\documentclass[11pt]{article}

\usepackage{graphicx}
\usepackage{epstopdf}
\usepackage{graphics}
\usepackage{amssymb}
\usepackage{amsmath}
\usepackage{amsthm}
\usepackage{amsfonts}
\usepackage{mathtools}
\usepackage{psfrag}
\usepackage{enumerate}
\usepackage{epsfig}
\usepackage{epsf}
\usepackage{latexsym}
\usepackage{fancyvrb}
\usepackage[small]{caption}
\usepackage{url}
\usepackage{array,multirow}
\usepackage{longtable}
\usepackage{float}
\usepackage{verbatim}
\usepackage{listings}
\usepackage{color}
\usepackage{tcolorbox}

\newtheorem{theorem}{Theorem}

\theoremstyle{remark}

\theoremstyle{plain}

\def\ie{{i.e.}}
\def\eg{{e.g.}}

\newcommand{\old}[1]{{}}
\newcommand{\later}[1]{{}}

\newcommand{\area}{{\rm area}}

\def\A{{\mathcal A}}

\def\L{{\mathcal L}}

\allowdisplaybreaks

\title{\textsc{New Lower Bounds for the Number of \\
Pseudoline Arrangements}}

\author{Adrian Dumitrescu\thanks{Department of Computer Science,
University of Wisconsin--Milwaukee, USA\@.
Email:~\texttt{dumitres@uwm.edu}.}
\and
Ritankar Mandal\thanks{Department of Computer Science,
University of Wisconsin--Milwaukee, USA\@.
Email:~\texttt{rmandal@uwm.edu}.}}

\begin{document}
\maketitle

\begin{abstract}
Arrangements of lines and pseudolines are fundamental objects 
in discrete and computational geometry. They also appear in 
other areas of computer science, such as the study of sorting
networks. Let $B_n$ be the number of nonisomorphic arrangements 
of $n$ pseudolines and let $b_n=\log_2{B_n}$. The problem of 
estimating $B_n$ was posed by Knuth in 1992. Knuth conjectured 
that $b_n \leq {n \choose 2} + o(n^2)$ and also derived the 
first upper and lower bounds: $b_n \leq 0.7924 (n^2 +n)$ and 
$b_n \geq n^2/6 -O(n)$. The upper bound underwent several 
improvements, $b_n \leq 0.6988\, n^2$ (Felsner, 1997), and 
$b_n \leq 0.6571\, n^2$ (Felsner and Valtr, 2011), for large 
$n$. Here we show that $b_n \geq cn^2 -O(n \log{n})$ for some 
constant $c>0.2083$. In particular, $b_n \geq 0.2083\, n^2$ 
for large $n$. This improves the previous best lower bound, 
$b_n \geq 0.1887\, n^2$, due to Felsner and Valtr (2011).
Our arguments are elementary and geometric in nature. Further, 
our constructions are likely to spur new developments and 
improved lower bounds for related problems, such as in
topological graph drawings. 

\medskip
\noindent
\textbf{\small Keywords}:
counting, pseudoline arrangement, recursive construction.
\end{abstract}

\section{Introduction} \label{sec:intro}

\begin{figure}[htpb]
\centering
\includegraphics[scale=0.51]{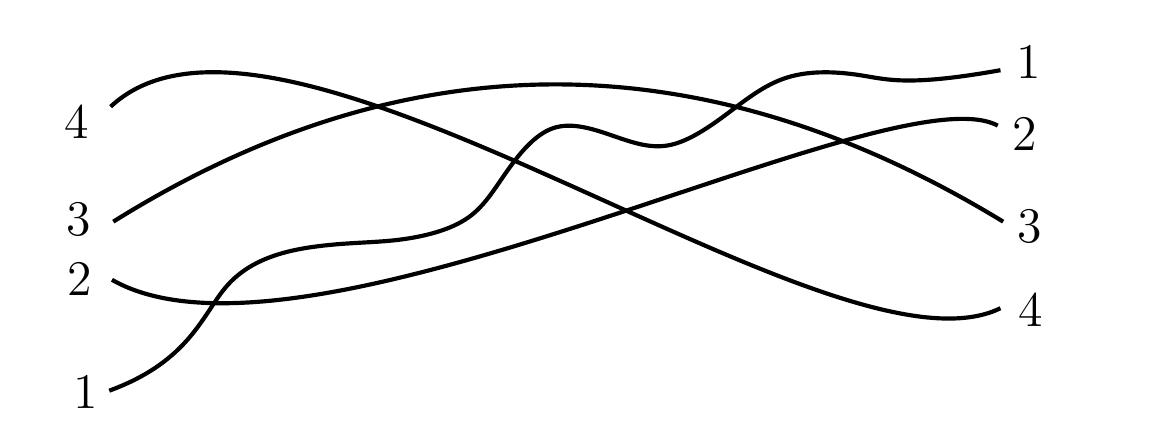}
\includegraphics[scale=0.51]{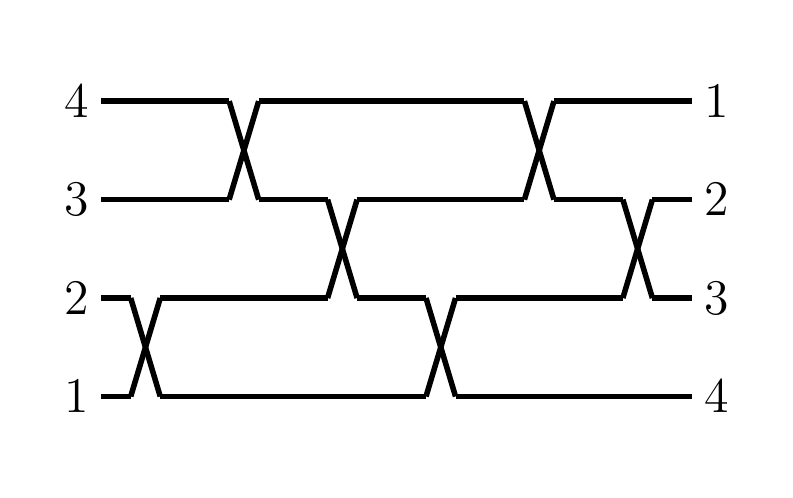}
\includegraphics[scale=0.51]{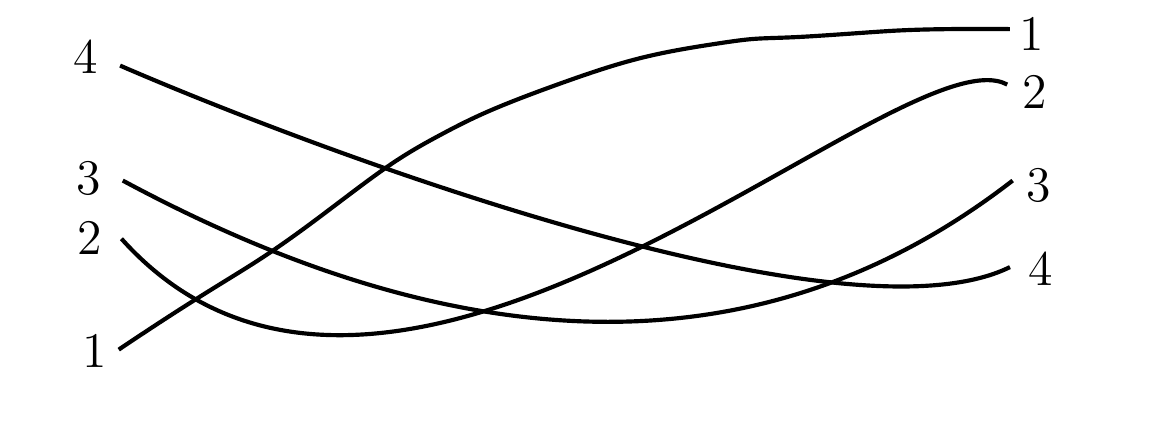}
\caption{Left: A simple arrangement $\A$. Center: Wiring 
diagram of $\A$. Right: An arrangement $\A'$ that is not 
isomorphic to the arrangement $\A$ on the left.}
\label{fig:r3r7r16}
\end{figure}

\paragraph{Arrangements of pseudolines.} A \emph{pseudoline} in 
the Euclidean plane is an $x$-monotone curve extending from 
negative infinity to positive infinity. An \emph{arrangement of 
pseudolines} is a family of pseudolines where each pair of 
pseudolines has a unique point of intersection (called 
`\emph{vertex}'). An arrangement is \emph{simple} if no three 
pseudolines have a common point of intersection, see 
Fig.~\ref{fig:r3r7r16}\,(left). Here the term arrangement always 
means simple arrangement if not specified otherwise. 

There are several representations and encodings of pseudoline 
arrangements. These representations help one count the number 
of arrangements. Three classic representations are \emph{allowable 
sequences} (introduced by Goodman and Pollack, see, 
\eg,~\cite{GP80,GP93}), 
\emph{wiring diagrams} (see for instance~\cite{FG17}), and 
\emph{zonotopal tilings} (see for instance~\cite{Fe04}). 
A~\emph{wiring diagram} is an Euclidean arrangement of pseudolines 
consisting of piece-wise linear `wires', each horizontal except 
for a short segment where it crosses another wire. Each pair of 
wires cross exactly once. The wiring diagram in 
Fig.~\ref{fig:r3r7r16}\,(center) 
represents the arrangement $\A$. The above representations have 
been shown to be equivalent; bijective proofs to this effect can 
be found in~\cite{Fe04}. Wiring diagrams are also known as 
\emph{reflection networks}, \ie, networks that bring $n$ wires 
labeled from $1$ to $n$ into their reflection by means of performing 
switches of adjacent wires; see~\cite[p.~35]{Kn92}. Lastly, they 
are also known under the name of \emph{primitive sorting networks};
see~\cite[Ch.~5.3.4]{Kn98}. The number of such networks  with $n$ wires,
\ie, the number of pseudoline arrangements with $n$ pseudolines,
is denoted by $A_n$. Stanley~\cite{St84} established the following
closed formula for $A_n$:
\[ A_n = \dfrac{\binom{n}{2} !}{\prod_{k=1}^{n-1} (2n - 2k - 1)^k}\]

Two arrangements are \emph{isomorphic} if they can be mapped onto 
each other by a homeomorphism of the plane; see Fig.~\ref{fig:r3r7r16}. 
The number of \emph{nonisomorphic} arrangements of $n$ pseudolines
is denoted by $B_n$; this is the number of equivalence classes
of all arrangements of $n$ pseudolines; see~\cite[p.~35]{Kn92}. 
This means that for $A_n$, the left to right order of the vertices
in the arrangement plays a role while for $B_n$ only the order 
of vertices along each particular pseudoline is important, \ie, 
the relative position of two vertices from distinct pairs of 
pseudolines does not matter. We are interested in the growth rate 
of $B_n$; so let\footnote{Throughout this paper, $\log{x}$ is the 
logarithm in base $2$ of $x$.} $b_n = \log_2 B_n$. Knuth~\cite{Kn92} 
conjectured that $b_n \leq {n \choose 2} + o(n^2)$; see also~\cite[p.~147]{FG17} 
and~\cite[p.~259]{Fe97}. This conjecture is still open. 

\paragraph{Upper bounds on the number of pseudoline arrangements.}
Felsner~\cite{Fe97} used a horizontal encoding of an arrangement 
in order to estimate $B_n$. An arrangement can be represented 
by a sequence of horizontal cuts. The $i$th cut is the list of 
the pseudolines crossing the $i$th pseudoline in the order of 
the crossings. Using this approach, Felsner~\cite[Thm.~1]{Fe97} 
obtained the upper bound $b_n \leq 0.7213(n^2 - n)$; he then
refined this bound by using \emph{replace matrices}. A replace 
matrix is a binary $n \times n$ matrix $M$ with the properties 
$\sum\nolimits_{j=1}^n m_{ij} = n - i$ for all $i$ and 
$m_{ij} \geq m_{ji}$ for all $i < j$. Using this technique, the author
established the upper bound $b_n < 0.6974\, n^2$~\cite[Thm.~2]{Fe97}.

In his seminal paper on the topic, Knuth~\cite{Kn92} took a 
vertical approach for encoding. Let $\A$ be an arrangement of 
$n$ pseudolines $\{\ell_1, \ldots, \ell_n\}$.
By adding pseudoline $\ell_{n+1}$ to $\A$, we get $\A'$, an arrangement 
of $n+1$ pseudolines. The course of $\ell_{n+1}$ describes 
a vertical cutpath from top to bottom. The number of cutpaths 
of $\A$ is exactly the number of arrangements $\A'$ such that 
$\A' \setminus \{\ell_{n+1}\}$ is isomorphic to $\A$. 
Let $\gamma_n$ denote the maximum number of cutpaths 
in an arrangement of $n$ pseudolines. Therefore, one has 
$B_{n+1} \leq \gamma_n \cdot B_n$; and $B_3=2$.
Knuth~\cite{Kn92} proved that $\gamma_n \leq 3^n$, concluding that
$B_n \leq 3^{\binom{n+1}{2}}$ and thus 
$b_n \leq 0.5\, (n^2+n) \log_2 3 \leq 0.7924\, (n^2 +n)$;
this computation can be streamlined so that it yields 
$b_n \leq 0.7924\, n^2$, see~\cite{FV11}. Knuth also conjectured 
that $\gamma_n \leq n \cdot 2^n$, but this was refuted by 
Ond\v{r}ej B\'ilka in 2010~\cite{FV11}; see also~\cite[p.~147]{FG17}.
Felsner and Valtr~\cite{FV11} proved a refined result,
$\gamma_n \leq 4n \cdot 2.48698^n$, by a careful analysis.
This  yields $b_n \leq 0.6571\, n^2$, which is the current best 
upper bound.

\paragraph{Lower bounds on the number of pseudoline arrangements.}
Knuth~\cite[p.~37]{Kn92} gave a recursive construction in the 
setting of reflection networks. The number of nonisomorphic 
arrangements in his construction, thus also $B_n$, obeys the recurrence
\[ B_n \geq 2^{n^2/8 -n/4} B_{n/2}. \]
By induction this yields $B_n \geq 2^{n^2/6 -5n/2}$.

\begin{figure}[htpb]
\centering
\includegraphics[scale=0.5]{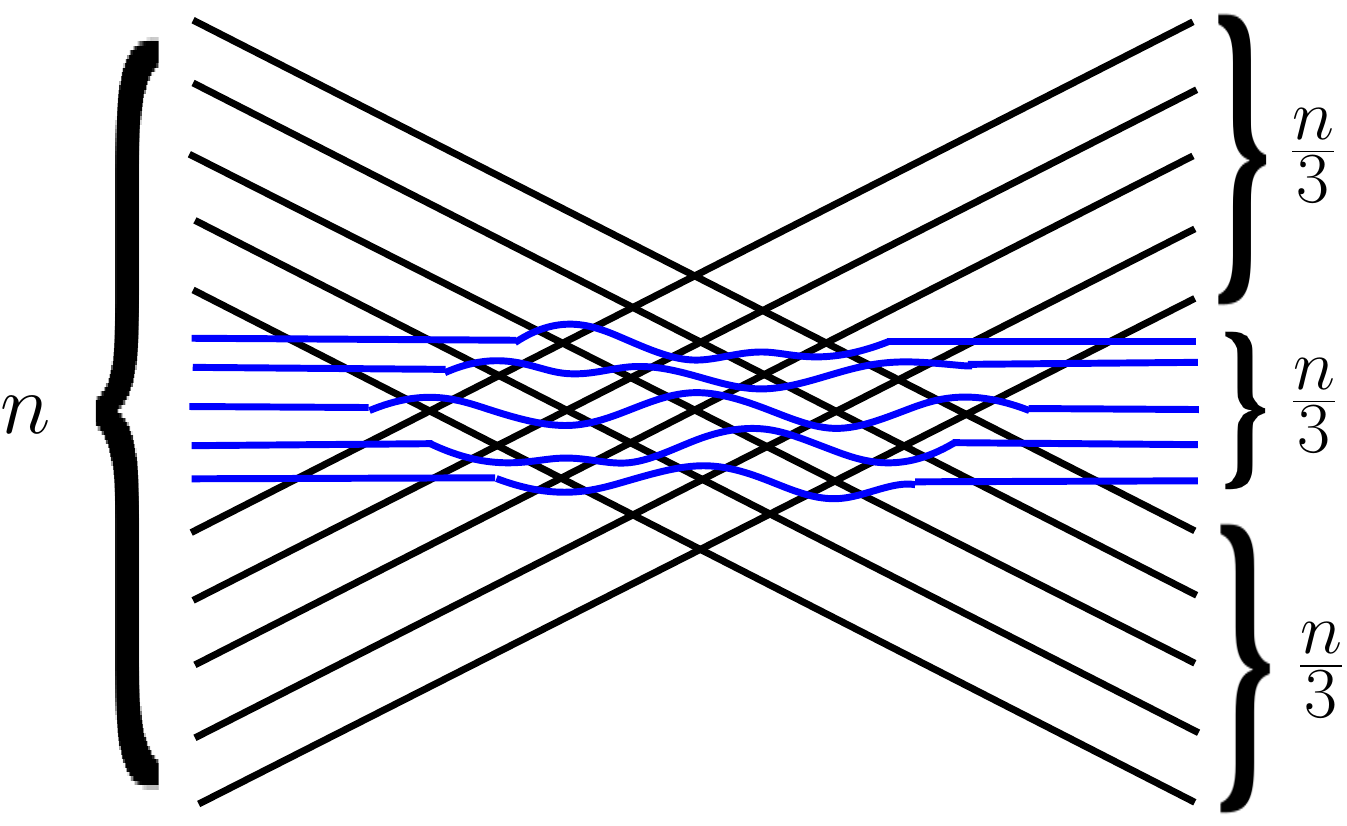}
\caption{Grid construction for a lower bound on $B_n$.}
\label{fig:matousek}
\end{figure}

Matou\v{s}ek sketched a simple---still recursive---grid construction 
in his book~\cite[Sec.~6.2]{Mat02}, see Fig~\ref{fig:matousek}. 
Let $n$ be a multiple of $3$ and $m = \frac{n}{3}$ (assume that 
$m$ is odd). The $2m$ lines in the two extreme bundles form a 
regular grid of $m^2$ points. The lines in the central bundle 
are incident to $\frac{3m^2+1}{4}$ of these grid points. At each 
such point, there are $2$ choices; going below it or above it, 
thus creating at least 
$\frac{3 m^2}{4} = \frac{3(n/3)^2}{4} = \frac{n^2}{12}$ binary choices.
Thus $B_n$ obeys the recurrence
\[ B_n \geq 2^{n^2/12} B_{n/3}^3, \]
which by induction yields $B_n \geq 2^{n^2/8}$.

\begin{figure}[htpb]
\centering
\includegraphics[scale=0.45]{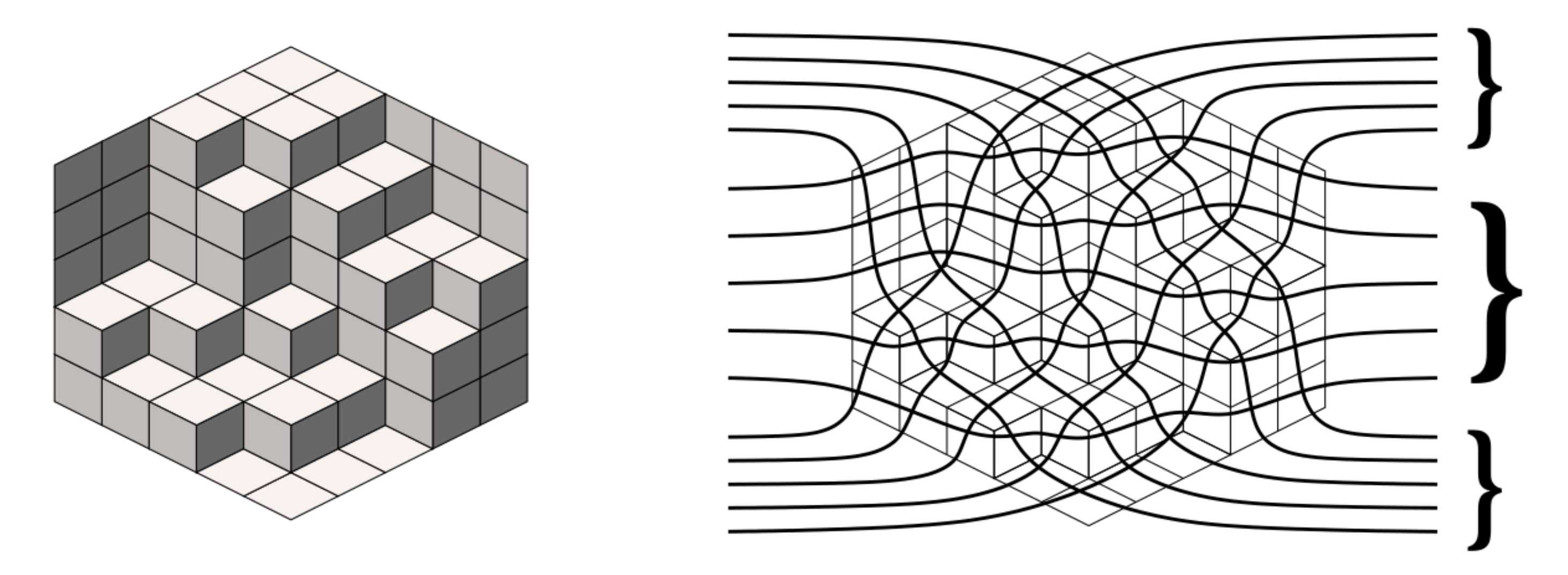}
\caption{The hexagon $H(5, 5, 5)$ with one of its rhombic tilings 
and a consistent partial arrangement corresponding to the tiling. 
This figure is reproduced from~\cite{FV11}.}
\label{fig:hex}
\end{figure}

Felsner and Valtr~\cite{FV11} used rhombic tilings of a centrally 
symmetric hexagon in an elegant recursive construction for a lower 
bound on $B_n$. Consider a set of $i+j+k$ pseudolines partitioned 
into the following three parts: 
$\{1, \ldots, i\}$, $\{i+1, \ldots, i+j\}$, $\{i+j+1, \ldots, i+j+k\}$, 
see Fig.~\ref{fig:hex}.
A partial arrangement is called \emph{consistent} if any two pseudolines 
from two different parts always cross but any two pseudolines from 
the same part never cross.

The zonotopal duals of consistent partial arrangements are rhombic 
tilings of the centrally symmetric hexagon $H(i,j,k)$ with side 
lengths $i,j,k$. The enumeration of rhombic tilings of $H(i, j, k)$ 
was solved by MacMahon~\cite{Mah60} (see also~\cite{Els84}), who 
showed that the number of tilings is
\begin{equation}\label{eq:macmahon}
P(i,j,k) = \prod\limits_{a=0}^{i-1} \prod\limits_{b=0}^{j-1}
\prod\limits_{c=0}^{k-1} \dfrac{a+b+c+2}{a+b+c+1}.
\end{equation}

A nontrivial (and quite involved) derivation using integral calculus 
shows that
\[ \ln P(n,n,n) = \left( \frac{9}{2} \ln 3 - 6 \ln 2 \right) n^2 + O(n \log{n})
= 1.1323\ldots n^2 +  O(n \log{n}). \] 
Assuming $n$ to be a multiple of $3$ in the recursion step, the 
construction yields the lower bound recurrence
\begin{equation}\label{eq:recurrencehex}
B_{n} \geq P \bigg(\frac{n}{3}, \frac{n}{3}, \frac{n}{3} \bigg) B_{n/3}^3.
\end{equation}
By induction, the analytic solution to formula~\eqref{eq:macmahon} 
together with the recurrence~\eqref{eq:recurrencehex} yield the 
lower bound $B_n \geq 2^{cn^2 -O(n \log^2{n})}$, where 
$c=\frac{3}{4} \ln 3 - \ln 2=0.1887\ldots$ In particular, 
$b_n \geq 0.1887\, n^2$ for large $n$; this is the previous best 
lower bound.

Table~\ref{tab:A_nB_n} shows the exact values of $A_n$ and $B_n$, and 
their growth rate (up to four digits after the decimal point) with 
respect to $n$, for small values of $n$. The values of $B_n$ for 
$n = 1$ to $9$ are from~\cite[p.~35]{Kn92} and the values of $B_{10}$, 
$B_{11}$, and $B_{12}$ are from~\cite{Fe97},~\cite{YN10}, and~\cite{SP95}, 
respectively; the values of $B_{13}$, $B_{14}$, and $B_{15}$ have been 
added recently, see~\cite{Ka11,SP95}. Observe that $A_n$ grows much 
faster than~$B_n$.
\begin{table}[H]
\begin{center}
\begin{tabular}{|r|r|l|r|l|}
\hline
$n$ & $A_n$ & $\dfrac{\log_2 A_n}{n^2}$ & $B_n$ & $\dfrac{\log_2 B_n}{n^2}$ \\
\hline \hline
$1$ & $1$ & $0$ & $1$ & $0$ \\ 
\hline
$2$ & $1$ & $0$ & $1$ & $0$ \\ 
\hline
$3$ & $2$ & $0.1111$ & $2$ & $0.1111$ \\ 
\hline
$4$ & $16$ & $0.25$ & $8$ & $0.1875$ \\ 
\hline
$5$ & $768$ & $0.3833$ & $62$ & $0.2381$ \\ 
\hline
$6$ & $292,864$ & $0.5044$ & $908$ & $0.2729$ \\ 
\hline
$7$ & $1,100,742,656$ & $0.6129$ & $24,698$ & $0.2977$ \\ 
\hline
$8$ & $48,608,795,688,960$ & $0.7104$ & $1,232,944$ & $0.3161$ \\ 
\hline
$9$ & $29,258,366,996,258,488,320$  & $0.7983$ & $112,018,190$ & $0.3301$ \\ 
\hline
$10$ &  &  & $18,410,581,880$ & $0.3409$ \\ 
\hline
$11$ &  &  & $5,449,192,389,984$ & $0.3496$ \\ 
\hline
$12$ &  &  & $2,894,710,651,370,536$ & $0.3566$ \\ 
\hline
$13$ &  &  & $2,752,596,959,306,389,652$ & $0.3624$ \\ 
\hline
$14$ &  &  & $4,675,651,520,558,571,537,540$ & $0.3672$ \\ 
\hline
$15$ &  &  & $14,163,808,995,580,022,218,786,390$ & $0.3713$ \\ 
\hline
\end{tabular}
\caption{Values of $A_n$ and $B_n$ for small $n$.}
\label{tab:A_nB_n}
\end{center}
\end{table}

\paragraph{Our results.}
Here we extend the method used by Matou\v{s}ek in his grid 
construction; observe that it uses lines of $3$ slopes. In 
Sections~\ref{sec:hexagonal6} (the 2nd part) and~\ref{sec:hexagonal12}, 
we use lines of $6$ and $12$ different slopes in hexagonal 
type constructions; yielding lower bounds $b_n \geq 0.1981\, n^2$ 
and $b_n \geq 0.2083\, n^2$ for large $n$, respectively.
In Sections~\ref{sec:rectangular8} and~\ref{sec:rectangular12} 
of the Appendix, we use lines of $8$ and $12$ different slopes 
in rectangular type constructions; yielding the lower bounds 
$b_n \geq 0.1999\, n^2$ and $b_n \geq 0.2053\, n^2$ for large 
$n$, respectively. 
While the construction in Section~\ref{sec:hexagonal12} gives 
a better bound, the one in Section~\ref{sec:hexagonal6} is 
easier to analyze. Results in Section~\ref{sec:hexagonal6},~\ref{sec:rectangular8} 
and~\ref{sec:rectangular12} are to appear in~\cite{DM19}. 
For each of the two styles, rectangular and hexagonal,
the constructions are presented in increasing order of complexity.
Our main result is summarized in the following.
\begin{theorem}\label{thm:main}
Let $B_n$ be the number of nonisomorphic arrangements of $n$ 
pseudolines. Then $B_n \geq 2^{cn^2 - O(n \log{n})}$, for some 
constant $c>0.2083$. In particular, $B_n \geq 2^{0.2083\, n^2}$ 
for large $n$.
\end{theorem}

\paragraph{Outline of the proof.}
We construct a line arrangement using lines of $k$ different slopes 
(for a small $k$). The final construction will be obtained by a 
small clockwise rotation, so that there are no vertical lines. 
Let $m=\lfloor n/k \rfloor$ or $m=\lfloor n/k \rfloor - 1$ 
(whichever is odd). Each bundle consists of $m$ equidistant lines 
in the corresponding parallel strip; remaining lines are discarded, 
or not used in the counting. An $i$-\emph{wise crossing}
is an intersection point of exactly $i$ lines. Let $\lambda_i(m)$ denote 
the number of $i$-wise crossings in the arrangement where each bundle 
consists of $m$ lines. Our goal is to estimate $\lambda_i(m)$ for 
each $i$. Then we can locally replace the lines around each $i$-wise 
crossing with any of the $B_i$ nonisomorphic pseudoline arrangements; 
and further apply recursively this construction to each of the $k$ 
bundles of parallel lines exiting this junction. This yields a simple 
pseudoline arrangement for each possible replacement choice. 
Consequently, the number of nonisomorphic pseudoline arrangements 
in this construction, say, $T(n)$, satisfies the recurrence:
\begin{equation} \label{eq:recurrence-k}
T(n) \geq F(n) \left[ T\left(\frac{n}{k}\right) \right]^k,
\end{equation} 
where $F(n)$ is a multiplicative factor counting the number of 
choices in this junction:
\begin{equation} \label{eq:f_n} 
F(n) \geq \prod_{i=3}^k B_i^{\lambda_i(n)}.
\end{equation} 

\paragraph{Related work.}
In a comprehensive recent paper, Kyn\v{c}l~\cite{Ky13} obtained 
estimates on the number of isomorphism classes of simple 
topological graphs that realize various graphs. The author 
remarks that it is probably hard to obtain tight estimates on this 
quantity, ``given that even for pseudoline arrangements, the best 
known lower and upper bounds on their number differ significantly''. 
While our improvements aren't spectacular, it seems however likely 
that some of the techniques we used here can be employed to obtain 
sharper lower bounds for topological graph drawings too.

\paragraph{Notations and formulas used.} For two similar figures 
$F,F'$, let $\rho(F,F')$ denote their similarity ratio. For a planar 
region $R$, let $\area(R)$ denote its area. By slightly abusing 
notation, let $\area(i,j,k)$ denote the area of the triangle made
by three lines $\ell_i$, $\ell_j$ and $\ell_k$. Assume that the 
equations of the three lines are $\alpha_s x + \beta_s y + \gamma_s = 0$, 
for $s=1,2,3$, respectively. Then
\begin{align*}
\area(i,j,k) &= \frac{A^2}{2 |C_1 C_2 C_3|}, \text { where } \\
A &= \begin{vmatrix}
\alpha_1 & \beta_1 & \gamma_1 \\ 
\alpha_2 & \beta_2 & \gamma_2 \\ 
\alpha_3 & \beta_3 & \gamma_3 
\end{vmatrix}, \\
C_1 &= (\alpha_2 \beta_3 - \beta_2 \alpha_3), \\
C_2 &= -(\alpha_1 \beta_3 - \beta_1 \alpha_3), \\
C_3 &= (\alpha_1 \beta_2 - \beta_1 \alpha_2).
\end{align*}

Let $P(i,j,g,h)$ denote the parallelogram made by the pairs of 
parallel lines $\ell_i \parallel \ell_j$ and $\ell_g \parallel \ell_h$.

\section{Preliminary constructions} \label{sec:hexagonal6}

\begin{figure}[htpb]
\centering
\includegraphics[scale=0.35]{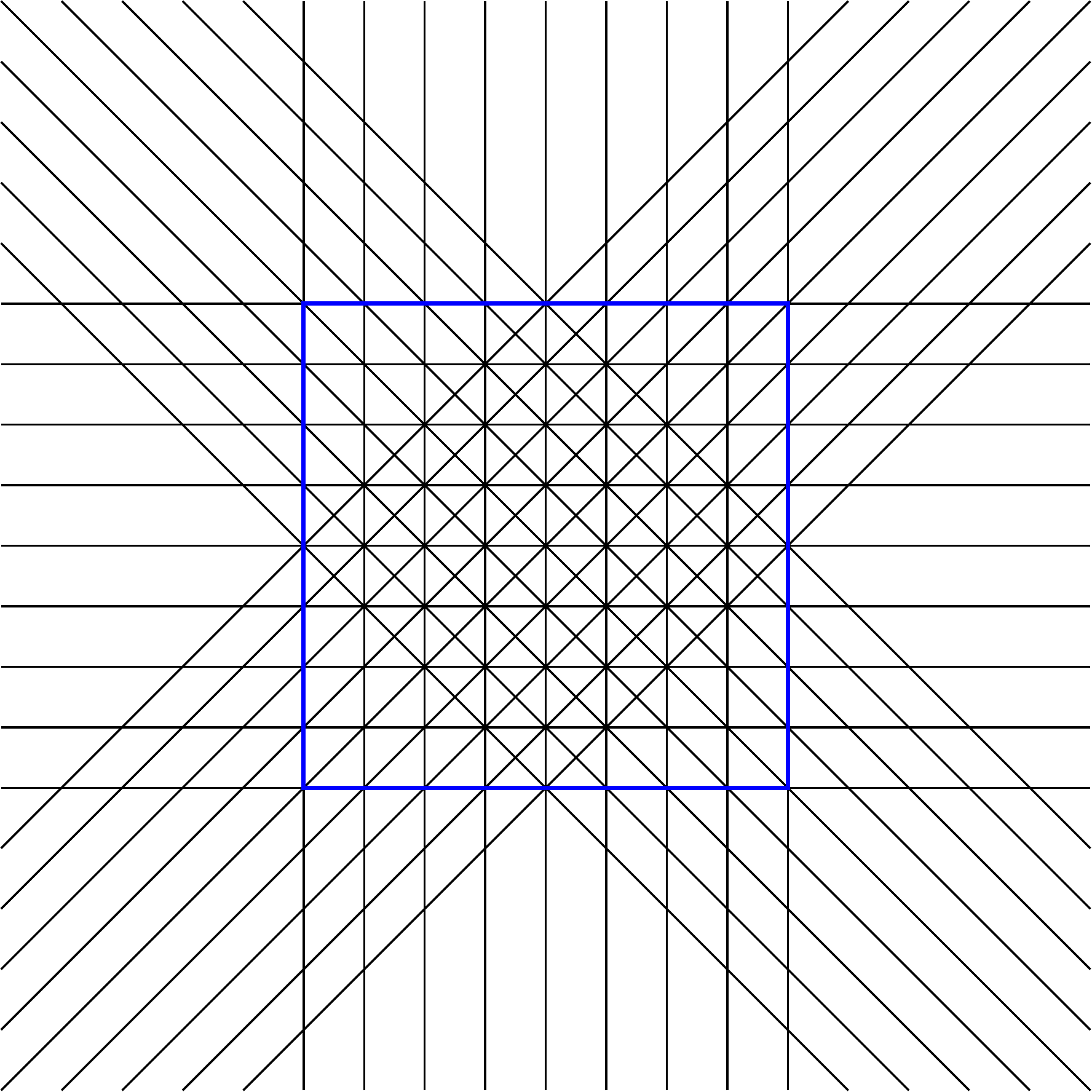}
\hspace*{1cm}
\includegraphics[scale=0.35]{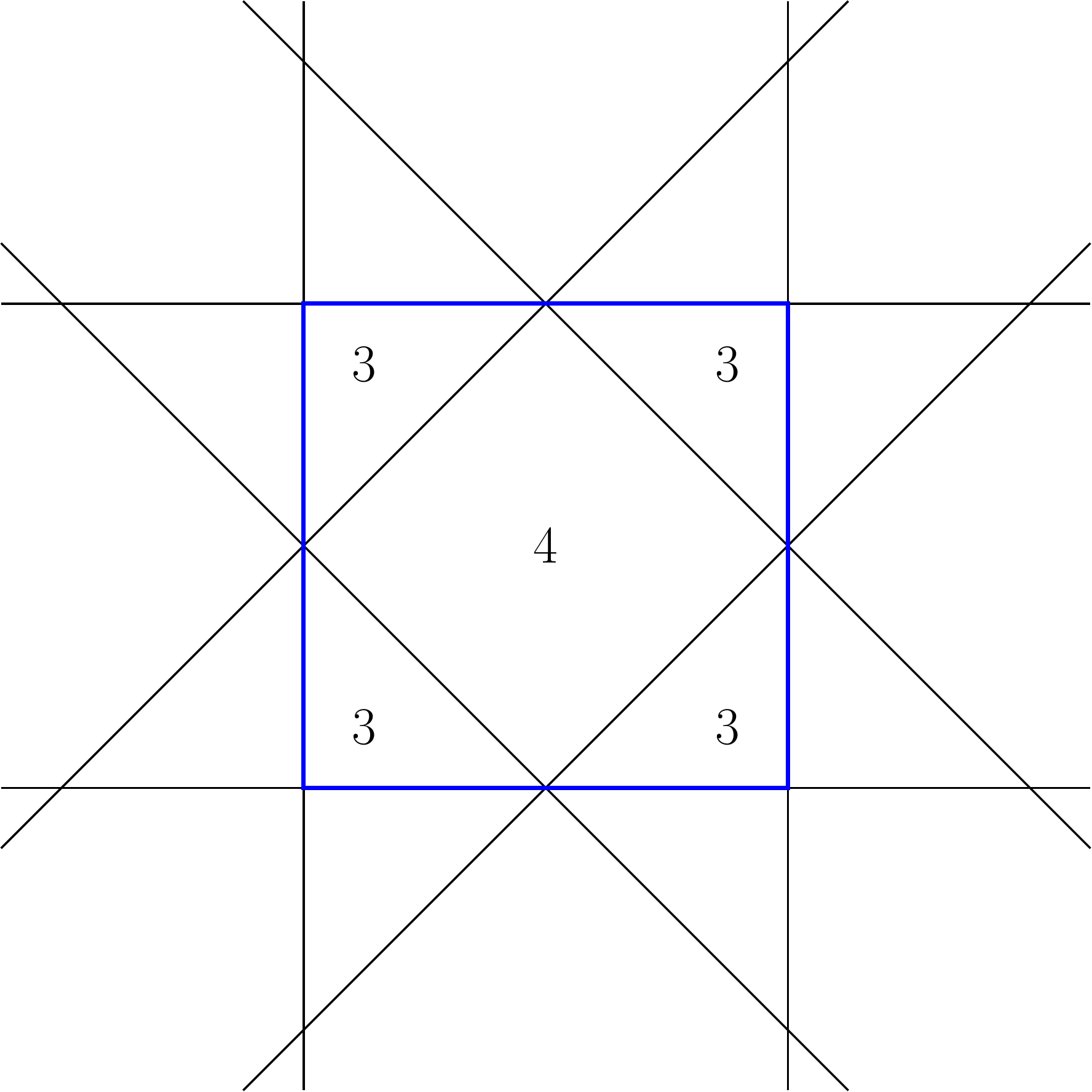}
\caption{Construction with $4$ slopes; here $m=9$. 
The unit square $U=[0,1]^2$ is shown in blue.}
\label{fig:f2}
\end{figure}

\paragraph{Warm-up: a rectangular construction with $4$ slopes.}
We start with a simple rectangular construction with $4$ bundles 
of parallel lines whose slopes are $0,\infty,\pm1$; see 
Fig.~\ref{fig:f2}. 
Let $U=[0,1]^2$ be the unit square we work with. The axes of all 
parallel strips are all incident to the center of $U$. 

For $i=3,4$, let $a_i$ denote the area of the region covered 
by exactly $i$ of the $4$ strips. It is easy to see that 
$a_3 = a_4 = 1/2$, and obviously $a_3 + a_4 = \area(U) = 1$. 
Observe that $\lambda_i(m)$ is proportional to $a_i$, for 
$i=3,4$; taking the boundary effect into account, we have
\[ \lambda_3(m) =  a_3 \, m^2 -O(m) = \frac{m^2}{2} -O(m), 
\text{ and \ \ } 
\lambda_4(m) =  a_4 \, m^2 -O(m) = \frac{m^2}{2} -O(m). \] 

Since $m = n/4$, $\lambda_i$ can be also viewed as a function of $n$.
Therefore $\lambda_3(n) = \lambda_4(n) = n^2/32 - O(n)$,
and so the multiplicative factor in Eq.\,\eqref{eq:recurrence-k}
is bounded from below as follows:  
\[ F(n) \geq \prod_{i=3}^4 B_i^{\lambda_i(n)} 
\geq 2^{n^2/32 - O(n)} \cdot 8^{n^2/32 - O(n)} 
= 2^{n^2/8 - O(n)}. \]
By induction on $n$, the resulting lower bound is 
$T(n) \geq 2^{n^2/6 - O(n \log{n})}$; this matches the 
constant $1/6$ in Knuth's lower bound described in 
Section~\ref{sec:intro}.

\begin{figure}[t]
\begin{minipage}[b]{0.65\linewidth}
\centering
\includegraphics[width=\linewidth]{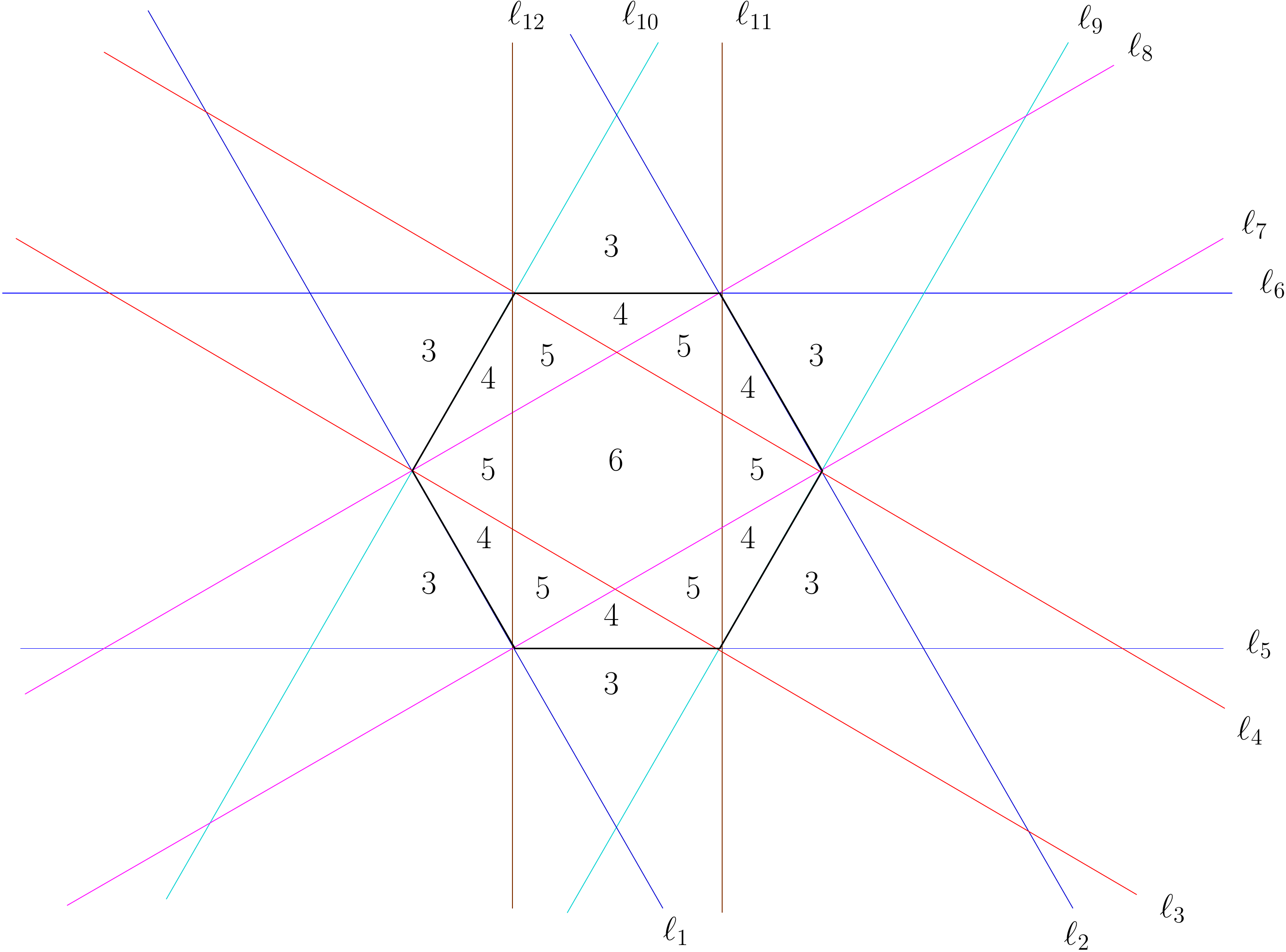}
\vspace{-3.7cm}
\end{minipage}
\hspace{0.01\linewidth}
\begin{minipage}[b]{0.26\linewidth}
\centering
\begin{tabular}{|r||c|c|c|} \hline
$i$ & $\alpha_i$ & $\beta_i$ & $\gamma_i$ \\ \hline \hline
$1$ & $\sqrt3$ & $1$ & $0$ \\ \hline
$2$ & $\sqrt3$ & $1$ & $-2\sqrt3$ \\ \hline
$3$ & $1$ & $\sqrt3$ & $-1$ \\ \hline 
$4$ & $1$ & $\sqrt3$ & $-3$ \\ \hline 
$5$ & $0$ & $1$ & $0$ \\ \hline 
$6$ & $0$ & $1$ & $-\sqrt3$ \\ \hline 
$7$ & $-1$ & $\sqrt3$ & $0$ \\ \hline 
$8$ & $-1$ & $\sqrt3$ & $-2$ \\ \hline 
$9$ & $-\sqrt3$ & $1$ & $\sqrt3$ \\ \hline 
$10$ & $-\sqrt3$ & $1$ & $-\sqrt3$ \\ \hline 
$11$ & $-1$ & $0$ & $1$ \\ \hline 
$12$ & $-1$ & $0$ & $0$ \\ \hline 
\end{tabular}
\end{minipage}
\caption{Left: The six parallel strips and corresponding 
covering multiplicities. These numbers only show incidences 
at the $3$-wise crossings made by primary lines. 
Right: Coefficients of the lines $\ell_i$ for $i=1,\ldots,12$.}
\label{fig:f9}
\end{figure}

\paragraph{Hexagonal construction with $6$ slopes.} 
We next describe and analyze a hexagonal construction with 
lines of $6$ slopes. Consider $6$ bundles of parallel lines 
whose slopes are $0,\infty, \pm1/\sqrt3, \pm\sqrt3$. Let $H$ 
be a regular hexagon whose side has unit length. Three parallel 
strips are bounded by the pairs of lines supporting opposite 
sides of $H$, while the other three parallel strips are bounded 
by the pairs of lines supporting opposite short diagonals of 
$H$. The axes of all six parallel strips are incident to the 
center of the circle created by the vertices of $H$; see 
Fig.~\ref{fig:f9}\,(left). This construction yields the lower 
bound $b_n \geq 0.1981\, n^2$ for large $n$. 

Assume a coordinate system where the lower left corner of $H$ 
is at the origin, and the lower side of $H$ lies along the $x$-axis.
Let $\L=\L_1 \cup \ldots \cup \L_6$ be the partition of $\L$ 
into six bundles of parallel lines. The $m$ lines in $\L_i$ 
are contained in the parallel strip bounded by the two lines 
$\ell_{2i-1}$ and $\ell_{2i}$, for $i=1,\ldots,6$. The equation 
of line $\ell_i$ is $\alpha_i x +\beta_i y +\gamma_i=0$, with
$\alpha_i,\beta_i,\gamma_i$, for $i=1,\ldots,12$, given in 
Fig~\ref{fig:f9}\,(right).

We refer to lines in $\L_1 \cup \L_3 \cup \L_5$ as the \emph{primary} 
lines, and to lines in $\L_2 \cup \L_4 \cup \L_6$ as \emph{secondary} 
lines. Note that
\begin{itemize} \itemsep 1pt
\item the distance between consecutive lines in any of the bundles 
of primary lines is
$\frac{\sqrt3}{m} \left(1-O\left(\frac{1}{m}\right)\right)$;
\item the distance between consecutive lines in any of the bundles 
of secondary lines is
$\frac{1}{m} \left(1-O\left(\frac{1}{m}\right)\right)$.
\end{itemize}

Let $\sigma_0=\sigma_0(m)$ and $\delta_0=\delta_0(m)$ denote the 
basic parallelogram and triangle respectively, determined by the 
consecutive lines in $\L_1 \cup \L_3 \cup \L_5$; the side length 
of $\sigma_0$ and $\delta_0$ is 
$\frac{2}{m} \left(1-O\left(\frac{1}{m}\right)\right)$. Let $H'$ 
be the smaller regular hexagon bounded by the short diagonals of $H$; 
the similarity ratio $\rho(H',H)$ is equal to $\frac{1}{\sqrt3}$.
Recall that (i)~the area of an equilateral triangle of side $s$ 
is $\frac{s^2 \sqrt3}{4}$; and (ii)~the area of a regular hexagon 
of side $s$ is $\frac{s^2 3\sqrt3}{2}$; as such, we have
\begin{align*}
\area(H) &= \frac{3 \sqrt3}{2}, \\
\area(H') &= \frac{\area(H)}{3} = \frac{\sqrt3}{2}, \\
\area(\delta_0) &= \frac{4}{m^2} \frac{\sqrt3}{4} \left(1-O\left(\frac{1}{m}\right)\right)
= \frac{\sqrt3}{m^2} \left(1-O\left(\frac{1}{m}\right)\right), \\
\area(\sigma_0) &= 2 \cdot \area(\delta_0)
= \frac{2\sqrt3}{m^2} \left(1-O\left(\frac{1}{m}\right)\right).
\end{align*}

For $i=3,4,5,6$, let $a_i$ denote the area of the (not necessarily 
connected) region covered by exactly $i$ of the $6$ strips. The 
following observations are in order: (i)~the six isosceles triangles 
based on the sides of $H$ inside $H$ have unit base and height 
$\frac{1}{2 \sqrt3}$; (ii)~the six smaller equilateral triangles 
incident to the vertices of $H$ have side-length $\frac{1}{\sqrt3}$. 
These observations yield
\begin{align*}
& a_3 = \area(H) = \frac{3 \sqrt3}{2}, \\
& a_4 = 6 \cdot \area(3,5,7) = 6 \cdot \frac{1}{4 \sqrt3} = \frac{\sqrt3}{2}, \\
& a_5 = 6 \cdot \area(3,7,11) = 6 \cdot \frac13 \frac{\sqrt3}{4} = \frac{\sqrt3}{2}, \\
& a_6 = \area(H') = \frac{\sqrt3}{2}. 
\end{align*}

Observe that $a_4 + a_5 + a_6 = \area(H)$. Recall that 
$\lambda_i(m)$ denote the number $i$-wise crossings where 
each bundle consists of $m$ lines. Note that $\lambda_i(m)$ 
is proportional to $a_i$, for $i=4,5,6$. Indeed, $\lambda_i(m)$ is 
equal to the number of $3$-wise crossings of lines in 
$\L_1 \cup \L_3 \cup \L_5$ that lie in a region covered by 
$i$ parallel strips, which is roughly equal to the ratio 
$\frac{a_i}{\area(\sigma_0)}$, for $i=4,5,6$. More precisely, 
taking also the boundary effect of the relevant regions into account, 
we obtain
\begin{align*}
& \lambda_4(m) = \frac{a_4}{\area(\sigma_0)} - O(m) 
 = \frac{\sqrt3}{2} \frac{m^2}{2 \sqrt3} - O(m) = \frac{m^2}{4} -O(m), \\
& \lambda_5(m) = \frac{a_5}{\area(\sigma_0)} - O(m) = \frac{m^2}{4} - O(m), \\
& \lambda_6(m) = \frac{a_6}{\area(\sigma_0)} - O(m) = \frac{m^2}{4} - O(m). 
\end{align*}

\begin{figure}[htpb]
\centering
\includegraphics[scale=0.4]{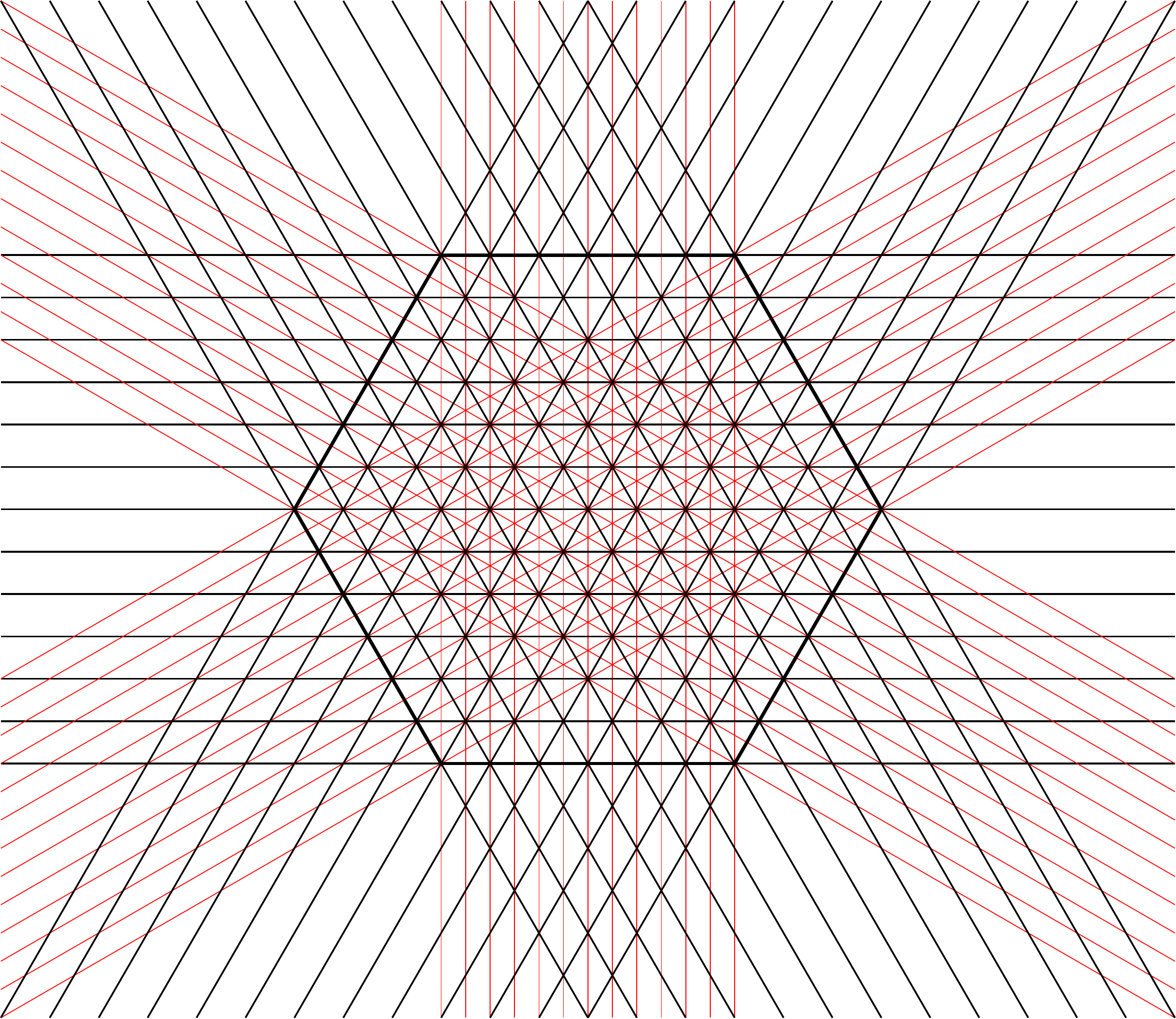}
\caption{Triple incidences of secondary lines (drawn in red).}
\label{fig:f12}
\end{figure}

For estimating $\lambda_3(m)$, the situation is little bit different, 
namely, in addition to considering $3$-wise crossings of the primary 
lines, we also observe $3$-wise crossings of the secondary lines at 
the centers of the small equilateral triangles contained in $H'$. 
See Fig.~\ref{fig:f12}. It follows that
\[ \lambda_3(m) = \frac{a_3}{\area(\sigma_0)} + \frac{\area(H')}{\area(\delta_0)} - O(m) 
= \frac{3 m^2}{4} + \frac{m^2}{2} - O(m) = \frac{5 m^2}{4} - O(m). \]
The values of $\lambda_i(m)$, for $i=3,4,5,6$, are summarized 
in Table~\ref{tab:lambda-6}; for convenience the linear terms are 
omitted. Since $m = n/6$, $\lambda_i$ can be also viewed as a function of $n$.
\renewcommand*{\arraystretch}{1.2}
\begin{table}[H]
\begin{center}
\begin{tabular}{|c||c|c|c|c|}
\hline
$i$ & $3$ & $4$ & $5$ & $6$ \\ \hline 
$\lambda_i(m)$ & $\frac{5 m^2}{4}$ & $\frac{m^2}{4}$ & $\frac{m^2}{4}$ 
& $\frac{m^2}{4}$ \\ \hline 
$\lambda_i(n)$ & $\frac{5 n^2}{4 \cdot 36}$ & $\frac{n^2}{4 \cdot 36}$
& $\frac{n^2}{4 \cdot 36}$ & $\frac{n^2}{4 \cdot 36}$ \\ \hline
\end{tabular}
\caption{The asymptotic values of $\lambda_i(m)$ and $\lambda_i(n)$ for $i=3,4,5,6$.}
\label{tab:lambda-6}
\end{center}
\end{table}
\renewcommand*{\arraystretch}{1}
The multiplicative factor in Eq.\,\eqref{eq:recurrence-k} is bounded 
from below as follows:  
\[ 
F(n) \geq \prod_{i=3}^6 B_i^{\lambda_i(n)} 
\geq 2^{5 n^2/144} \cdot 8^{n^2/144} \cdot 62^{n^2/144} \cdot 908^{n^2/144} \cdot 2^{-O(n)}. 
\]

We prove by induction on $n$ that $T(n) \geq 2^{cn^2 - O(n \log{n})}$ 
for a suitable constant $c>0$. It suffices to choose $c$ (using the 
values of $B_i$ for $i=3,4,5,6$ in Table~\ref{tab:A_nB_n}) so that 
\[ \frac{8 + \log{62} + \log{908}}{144} \geq \frac{5c}{6}. \]
The above inequality holds if we set
$c = \dfrac{\log(256 \cdot 62 \cdot 908)}{120} > 0.1981$,
and the lower bound follows.

\section{Hexagonal construction with $12$ slopes} \label{sec:hexagonal12}

\begin{figure}[ht]
\centering
\includegraphics[scale=0.9]{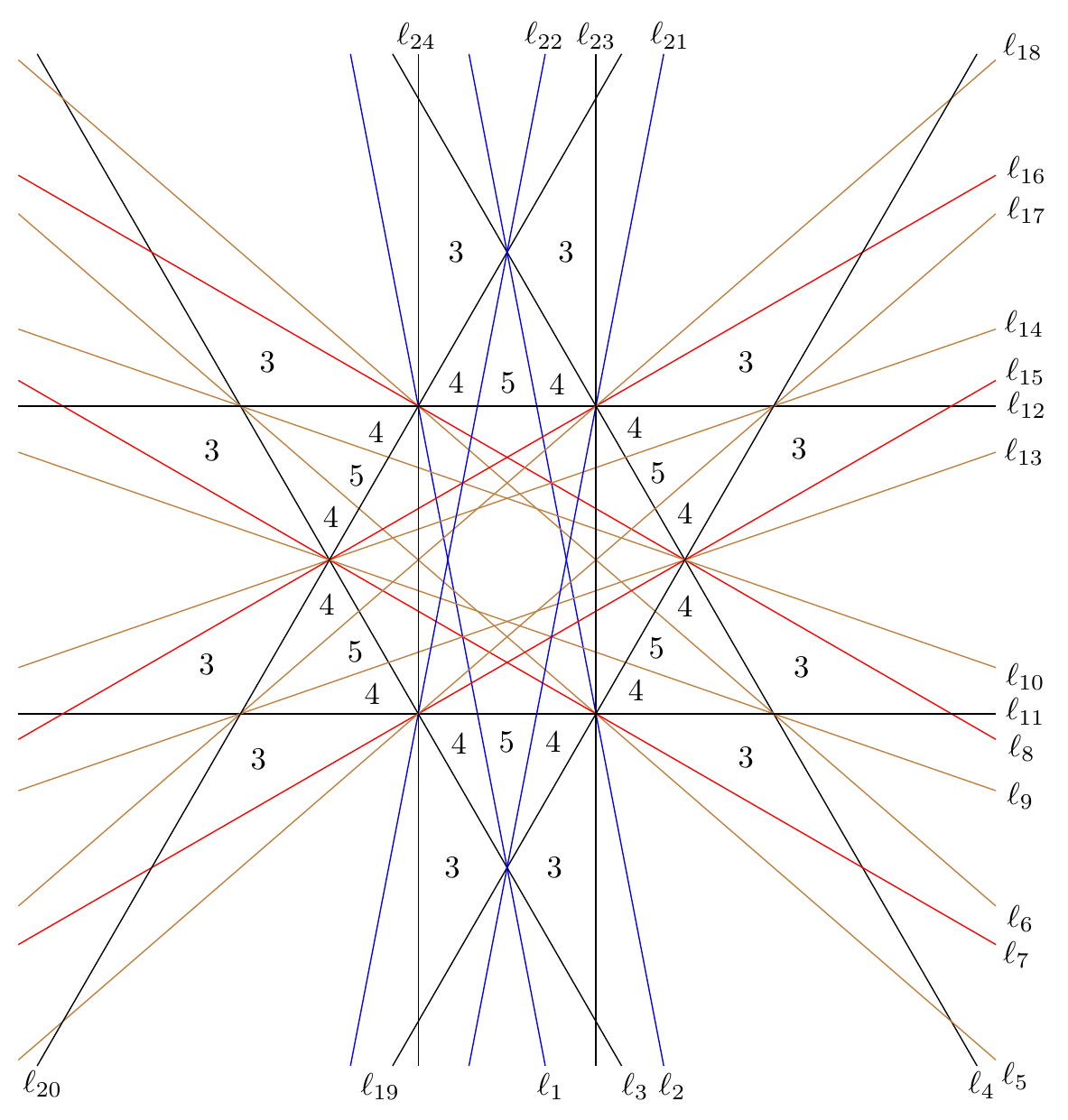}
\caption{Construction with $12$ slopes. The twelve parallel strips 
and the corresponding covering multiplicities. These numbers 
only reflect incidences at the grid vertices made by the primary lines. 
The numbers inside $H$ are shown in Fig.~\ref{fig:f16}}
\label{fig:f15}
\end{figure}

\begin{figure}[ht]
\centering
\includegraphics[scale=0.78]{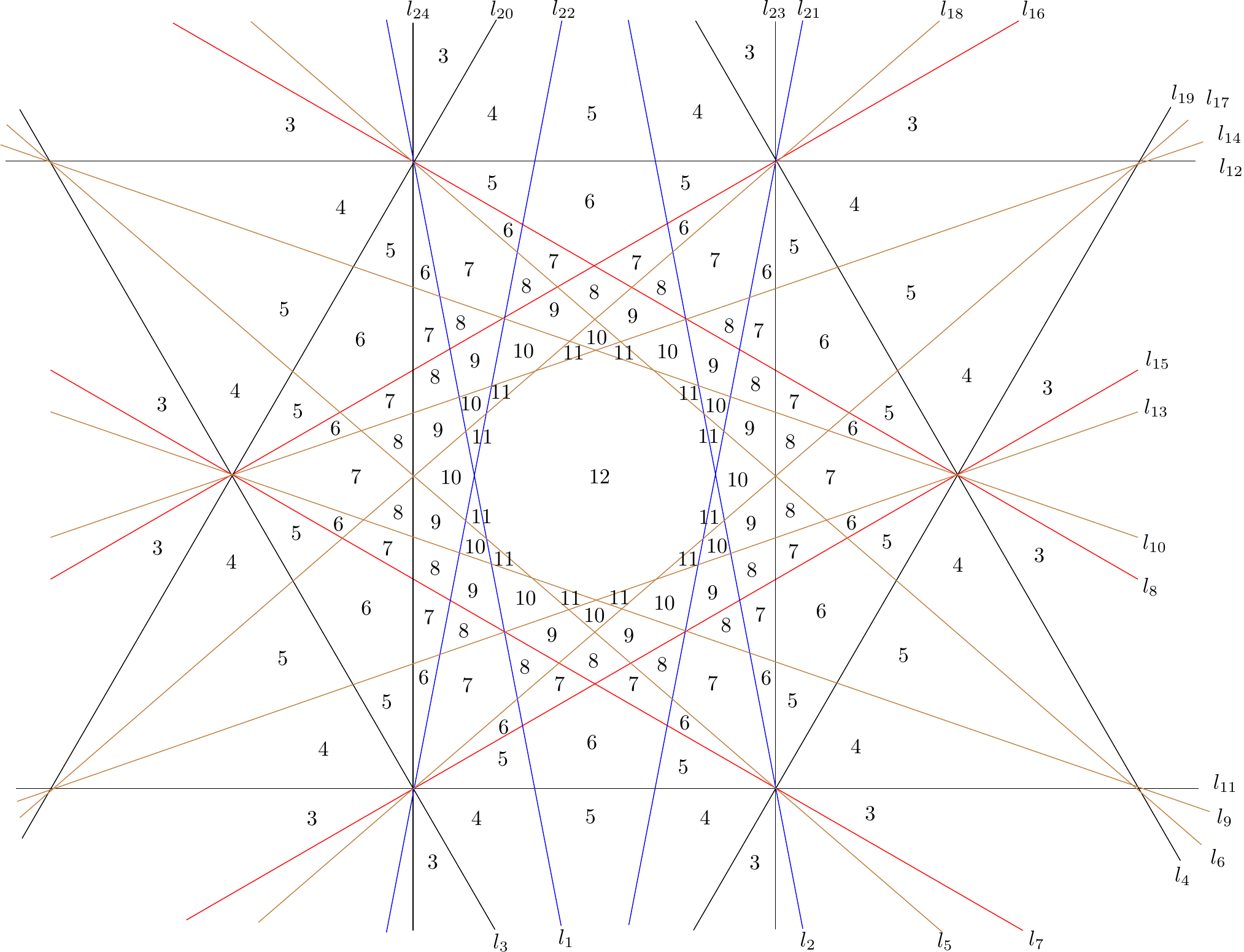}
\caption{Construction with $12$ slopes. The twelve parallel 
strips and the corresponding covering multiplicities. These 
numbers only reflect incidences at the grid vertices made 
by the primary lines.}
\label{fig:f16}
\end{figure}

We next describe and analyze a hexagonal construction with 
lines of $12$ slopes, which provides our main result in 
Theorem~\ref{thm:main}. Consider $12$ bundles of parallel 
lines whose slopes are $0, \infty, \pm\sqrt3/5$, $\pm1/\sqrt3, 
\pm\sqrt3/2, \pm\sqrt3, \pm3\sqrt3$. Let $H$ be a regular 
hexagon whose side has unit length. The axes of all parallel 
strips are incident to the center of the circle created by 
the vertices of $H$; see Figs.~\ref{fig:f15} and~\ref{fig:f16}. 
This construction yields the lower bound $b_n \geq 0.2083\, n^2$ 
for large $n$.

\renewcommand*{\arraystretch}{1.1}
\begin{table}[htpb]
\begin{center}
\begin{minipage}{0.32\linewidth}
\centering
\begin{tabular}{|r||c|c|c|} \hline
$i$ & $\alpha_i$ & $\beta_i$ & $\gamma_i$ \\ \hline \hline
$1$ & $3\sqrt3$ & $1$ & $-\sqrt3$ \\ \hline 
$2$ & $3\sqrt3$ & $1$ & $-3\sqrt3$ \\ \hline 
$3$ & $\sqrt{3}$ & $1$ & $0$ \\ \hline
$4$ & $\sqrt{3}$ & $1$ & $-2\sqrt{3}$ \\ \hline
$5$ & $\sqrt3$ & $2$ & $-\sqrt3$ \\ \hline 
$6$ & $\sqrt3$ & $2$ & $-2\sqrt3$ \\ \hline
$7$ & $1$ & $\sqrt{3}$ & $-1$ \\ \hline 
$8$ & $1$ & $\sqrt{3}$ & $-3$ \\ \hline
\end{tabular}
\end{minipage}
\begin{minipage}{0.32\linewidth}
\centering
\begin{tabular}{|r||c|c|c|} \hline
$i$ & $\alpha_i$ & $\beta_i$ & $\gamma_i$ \\ \hline \hline
$9$ & $\sqrt3$ & $5$ & $-2\sqrt3$ \\ \hline 
$10$ & $\sqrt3$ & $5$ & $-4\sqrt3$ \\ \hline
$11$ & $0$ & $1$ & $0$ \\ \hline 
$12$ & $0$ & $1$ & $-\sqrt{3}$ \\ \hline 
$13$ & $-\sqrt3$ & $5$ & $-\sqrt3$ \\ \hline 
$14$ & $-\sqrt3$ & $5$ & $-3\sqrt3$ \\ \hline
$15$ & $-1$ & $\sqrt3$ & $0$ \\ \hline 
$16$ & $-1$ & $\sqrt3$ & $-2$ \\ \hline
\end{tabular}
\end{minipage}
\begin{minipage}{0.32\linewidth}
\centering
\begin{tabular}{|r||c|c|c|} \hline
$i$ & $\alpha_i$ & $\beta_i$ & $\gamma_i$ \\ \hline \hline 
$17$ & $-\sqrt3$ & $2$ & $0$ \\ \hline 
$18$ & $-\sqrt3$ & $2$ & $-\sqrt3$ \\ \hline
$19$ & $-\sqrt3$ & $1$ & $\sqrt3$ \\ \hline 
$20$ & $-\sqrt3$ & $1$ & $-\sqrt3$ \\ \hline 
$21$ & $-3\sqrt3$ & $1$ & $2\sqrt3$ \\ \hline 
$22$ & $-3\sqrt3$ & $1$ & $0$ \\ \hline 
$23$ & $-1$ & $0$ & $1$ \\ \hline 
$24$ & $-1$ & $0$ & $0$ \\ \hline 
\end{tabular}
\end{minipage}
\end{center}
\caption{Coefficients of the $24$ lines.}
\label{tab:24lines}
\end{table}
\renewcommand*{\arraystretch}{1}

Assume a coordinate system where the lower left corner of 
$H$ is at the origin, and the lower side of $H$ lies along 
the $x$-axis. Let $\L=\L_1 \cup \ldots \cup \L_{12}$ be the 
partition of $\L$ into twelve bundles of parallel lines. 
The $m$ lines in $\L_i$ are contained in the parallel strip 
$\Gamma_i$ bounded by the two lines $\ell_{2i-1}$ and 
$\ell_{2i}$, for $i=1,\ldots,24$. The equation of line 
$\ell_i$ is $\alpha_i x +\beta_i y +\gamma_i=0$, with
$\alpha_i,\beta_i,\gamma_i$, for $i=1,\ldots,24$, given 
in Table~\ref{tab:24lines}. 

$\Gamma_2$, $\Gamma_6$ and $\Gamma_{10}$ are bounded by the pairs of 
lines supporting opposite sides of $H$, while 
$\Gamma_4$, $\Gamma_8$ and  $\Gamma_{12}$ are bounded by the pairs of 
lines supporting opposite short diagonals of $H$. Therefore 
$H = \Gamma_2 \cap \Gamma_6 \cap \Gamma_{10}$. 
We refer to lines in $\L_2 \cup \L_6 \cup \L_{10}$ as the 
\emph{primary} lines, to lines in $\L_4 \cup \L_8 \cup \L_{12}$ 
as the \emph{secondary} lines, and to the rest of the lines 
as the \emph{tertiary} lines. Note that  

\begin{itemize} \itemsep 1pt
\item the distance between consecutive lines in any of the 
bundles of primary lines is
$\frac{\sqrt3}{m} \left(1-O\left(\frac{1}{m}\right)\right)$;
\item the distance between consecutive lines in any of the 
bundles of secondary lines is
$\frac{1}{m} \left(1-O\left(\frac{1}{m}\right)\right)$;
\item the distance between consecutive lines in any of the 
bundles of tertiary lines is
$\sqrt{\frac37}\frac{1}{m} \left(1-O\left(\frac{1}{m}\right)\right)$.
\end{itemize}

We refer to the intersection points of the primary lines as 
\emph{grid vertices}. There are two types of grid vertices:
the grid vertices in $H$ are intersection of $3$ primary lines 
and the ones outside $H$ are intersection of $2$ primary lines.

Let $\sigma_0=\sigma_0(m)$ and $\delta_0=\delta_0(m)$ denote 
the basic parallelogram and triangle respectively, determined 
by the primary lines (\ie, lines in $\L_2 \cup \L_6 \cup \L_{10}$); 
the side length of $\sigma_0(m)$ and $\delta_0$ is 
$\frac{2}{m} \left(1-O\left(\frac{1}{m}\right)\right)$. 
We refer to these basic parallelograms as \emph{grid cells}.
Recall that (i)~the area of an equilateral triangle of side 
$s$ is $\frac{s^2 \sqrt3}{4}$; and (ii)~the area of a regular 
hexagon of side $s$ is $\frac{s^2 3\sqrt3}{2}$; as such, we have
\begin{align*}
\area(H) &= \frac{3 \sqrt3}{2}, \\
\area(\delta_0) &= \frac{4}{m^2} \frac{\sqrt3}{4} \left(1-O\left(\frac{1}{m}\right)\right)
= \frac{\sqrt3}{m^2} \left(1-O\left(\frac{1}{m}\right)\right), \\
\area(\sigma_0) &= 2 \cdot \area(\delta_0) =
\frac{2\sqrt3}{m^2} \left(1-O\left(\frac{1}{m}\right)\right).
\end{align*}

For $i=3,\ldots,12$, let $a_i$ denote the area of the (not 
necessarily connected) region covered by exactly $i$ of 
the $6$ strips. Recall that $\area(i,j,k)$ denotes the area 
of the triangle made by $\ell_i$, $\ell_j$ and $\ell_k$. 

Observe that $a_{12}$ is the area of the $12$-gon $\bigcap_{i=1}^{12} \Gamma_i$. 
This $12$-gon is not regular, since consecutive vertices lie on 
two concentric cycles of radii $\frac{1}{3}$ and $\frac{\sqrt3}{5}$ 
centered at $(\frac{1}{2}, \frac{\sqrt3}{2})$. So $a_{12}$ is the 
sum of the areas of $12$ congruent triangles; each with one vertex 
at the center of $H$ and other two as the two consecutive vertices 
of the $12$-gon. Each of these triangles have area $\frac{\sqrt3}{60}$. 
Therefore, 
\begin{align*}
a_{12} & = 12 \cdot \frac{\sqrt3}{60} = \frac{\sqrt3}{5}, \\
a_{11} & = 12 \cdot \area(1,5,9) = 12 \cdot \frac{1}{140\sqrt3} = \frac{\sqrt3}{35}, \\
a_{10} & = 6 \cdot ( \area(1,5,13) - \area(1,5,9) ) + 6 \cdot ( \area(5,9,22) - \area(1,5,9) ) \\
& = 6 \cdot \left( \frac{\sqrt3}{70} - \frac{1}{140\sqrt3} \right) 
+ 6 \cdot \left( \frac{1}{56\sqrt3} - \frac{1}{140\sqrt3} \right) = \frac{13\sqrt3}{140}, \\
a_9 & = 12 \cdot ( \area(1,7,22) - \area(1,9,22) ) 
= 12 \cdot \left( \frac{1}{20\sqrt3} - \frac{1}{56\sqrt3} \right) = \frac{9\sqrt3}{70}, \\
a_8 & = 6 \cdot ( \area(9,22,24) - \area(7,22,24) ) + 12 \cdot \area(7,13,22) \\
& = 6 \cdot \left( \frac{\sqrt3}{40} - \frac{1}{20\sqrt3} \right)
+ 12 \cdot \frac{\sqrt3}{140} = \frac{19\sqrt3}{140}, \\
a_7 & = 12 \cdot ( \area(7,22,24) - \area(13,22,24) ) 
+ 6 \cdot ( \area(1,17,22) - \area(1,13,22) ) \\
& = 12 \cdot \left( \frac{1}{20\sqrt3} - \frac{\sqrt3}{140} \right) 
+ 6 \cdot \left( \frac{5}{28\sqrt3} - \frac{1}{14\sqrt3} \right) = \frac{23\sqrt3}{70}, \\
a_6 & = 12 \cdot ( \area(13,22,24) ) + 6 \cdot ( \area(7,11,15) - 2 \cdot \area(1,11,15) ) \\
& = 12 \cdot \frac{\sqrt3}{140} + 6 \cdot \left( \frac{1}{4\sqrt3}
- 2 \cdot \frac{1}{20\sqrt3} \right) = \frac{27\sqrt3}{70}, \\
a_5 & = 12 \cdot ( \area(1,11,15) ) + 6 \cdot ( \area(1,11,21) ) 
= 12 \cdot \frac{1}{20\sqrt3} + 6 \cdot \frac{1}{4\sqrt3} = \frac{7\sqrt3}{10}, \\
a_4 & = 12 \cdot ( \area(1,3,11) ) = 12 \cdot \frac{1}{4\sqrt3} = \sqrt3, \\
a_3 & = 12 \cdot ( \area(4,7,11) ) = 12 \cdot \frac{\sqrt3}{4} = 3\sqrt3. 
\end{align*}

Observe that the region whose area is $\sum_{i=5}^{12} a_i$ consists of
the hexagon $H$ and $6$ triangles outside $H$. Therefore, 
\[ \sum\nolimits_{i=5}^{12} a_i = \area(H) + 6 \cdot \area(1,11,21) 
= \frac{3 \sqrt3}{2} + 6 \cdot \frac{1}{4\sqrt3} = 2\sqrt3. \]

Recall that $\lambda_i(m)$ denotes the number of $i$-wise 
crossings where each bundle consists of $m$ lines. Note that 
$\lambda_i(m)$ is proportional to $a_i$, for $i=5,6,\ldots,12$. 
Indeed, $\lambda_i(m)$ is equal to the number of grid vertices 
that lie in a region covered by $i$ parallel strips, which is 
roughly equal to the ratio $\frac{a_i}{\area(\sigma_0)}$, for 
$i=5,6,\ldots,12$. More precisely, taking also the boundary 
effect of the relevant regions into account, we obtain
\begin{align*}
\lambda_{12}(m) &= \frac{a_{12}}{\area(\sigma_0)} - O(m) 
= \frac{\sqrt3}{5} \frac{m^2}{2\sqrt3} - O(m) = \frac{m^2}{10} - O(m), \\
\lambda_{11}(m) &= \frac{a_{11}}{\area(\sigma_0)} - O(m) = \frac{m^2}{70} - O(m), \\
\lambda_{10}(m) &= \frac{a_{10}}{\area(\sigma_0)} - O(m) = \frac{13 m^2}{280} - O(m), \\
\lambda_9(m) &= \frac{a_9}{\area(\sigma_0)} - O(m) = \frac{9 m^2}{140} - O(m), \\
\lambda_8(m) &= \frac{a_8}{\area(\sigma_0)} - O(m) = \frac{19 m^2}{280} - O(m), \\
\lambda_7(m) &= \frac{a_7}{\area(\sigma_0)} - O(m) = \frac{23 m^2}{140} - O(m), \\
\lambda_6(m) &= \frac{a_6}{\area(\sigma_0)} - O(m) = \frac{27 m^2}{140} - O(m), \\
\lambda_5(m) &= \frac{a_6}{\area(\sigma_0)} - O(m) = \frac{7 m^2}{20} - O(m).
\end{align*}

\begin{figure}[htpb]
\centering
\includegraphics[scale=0.82]{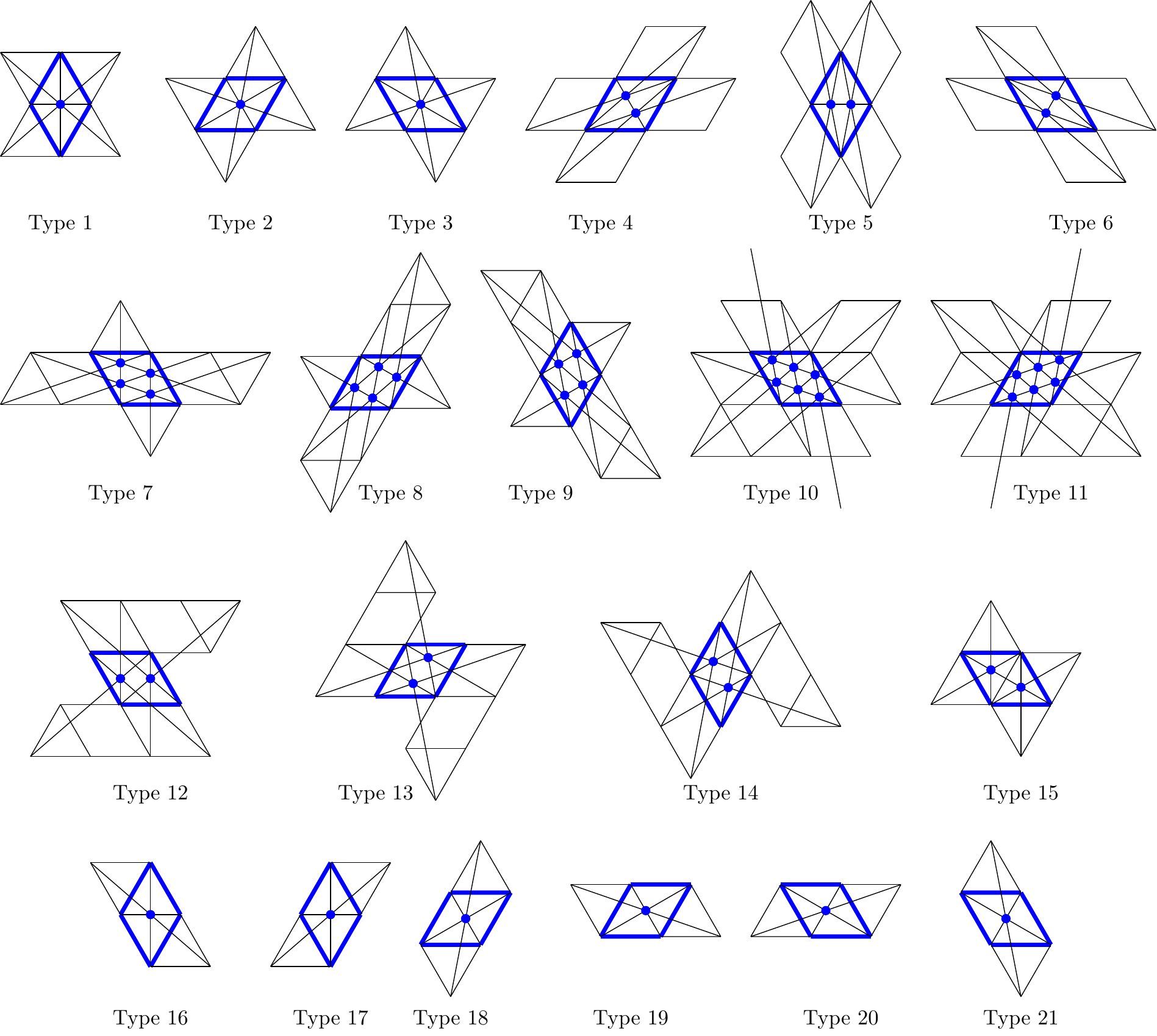}
\caption{Types of incidences of $3$ and $4$ lines that are not at grid vertices: 
$4$-wise crossings: types $1$ through $3$; 
$3$-wise crossings: types $4$ through $21$. 
To list the coordinates of the crossing points (shown as blue dots), 
we set the leftmost vertex of the grid cell (shown in blue lines) at $(0,0)$ 
and the length of the sides of each grid cell as $1$. \\
For types $1$ through $3$ the crossings are at the center of the parallelogram. \\ 
For types $4$ through $6$, the crossings are on the short diagonal at 
$\frac{1}{3}$rd and $\frac{2}{3}$rd of the short diagonal. \\ 
For type $7$, the crossings are at $(\frac{1}{2}, \frac{-\sqrt3}{10})$, 
$(\frac{1}{2}, \frac{-3\sqrt3}{10})$, $(1, \frac{-\sqrt3}{5})$, 
$(1, \frac{-2\sqrt3}{5})$. \\ 
For type $8$, the crossings are at $(\frac{2}{5}, \frac{\sqrt3}{5})$, 
$(\frac{7}{10}, \frac{\sqrt3}{10})$, $(\frac{4}{5}, \frac{2\sqrt3}{5})$, 
$(\frac{11}{10}, \frac{3\sqrt3}{10})$. \\
For type $9$, the crossings are at $(\frac{3}{10}, \frac{\sqrt3}{10})$, 
$(\frac{2}{5}, \frac{-\sqrt3}{5})$, $(\frac{3}{5}, \frac{\sqrt3}{5})$, 
$(\frac{7}{10}, \frac{-\sqrt3}{10})$. \\ 
For type $10$, the crossings are at $(\frac{5}{14}, \frac{-\sqrt3}{14})$, 
$(\frac{5}{7}, \frac{-\sqrt3}{7})$, $(\frac{15}{14}, \frac{-3\sqrt3}{14})$, 
$(\frac{3}{7}, \frac{-2\sqrt3}{7})$, $(\frac{11}{14}, \frac{-5\sqrt3}{14})$, 
$(\frac{8}{7}, \frac{-3\sqrt3}{7})$. \\ 
For type $11$, the crossings are at $(\frac{8}{7}, \frac{3\sqrt3}{7})$, 
$(\frac{11}{14}, \frac{5\sqrt3}{14})$, $(\frac{3}{7}, \frac{2\sqrt3}{7})$, 
$(\frac{15}{14}, \frac{3\sqrt3}{14})$, $(\frac{5}{7}, \frac{\sqrt3}{7})$, 
$(\frac{5}{14}, \frac{\sqrt3}{14})$. \\ 
For type $12$, the crossings are at $(\frac{1}{2}, \frac{-\sqrt3}{4})$ 
and $(1, \frac{-\sqrt3}{4})$. \\ 
For type $13$, the crossings are at $(\frac{5}{8}, \frac{\sqrt3}{8})$ 
and $(\frac{7}{8}, \frac{3\sqrt3}{8})$. \\ 
For type $14$, the crossings are at $(\frac{3}{8}, \frac{\sqrt3}{8})$ 
and $(\frac{5}{8}, \frac{-\sqrt3}{8})$. \\ 
For type $15$, the crossings are on the long diagonal at $\frac{1}{3}$rd 
and $\frac{2}{3}$rd of the long diagonal. \\
For types $16$ through $21$ the crossings are at the center of the parallelogram. \\  
The relative positions of all these crossings are shown in Fig.~\ref{fig:f18}.}
\label{fig:f17}
\end{figure}

For $i=3,4$, not all the $i$-wise crossings are at grid 
vertices. It can be exhaustively verified (by hand) that 
there are $21$ types of crossings; see Fig.~\ref{fig:f17}. 
Types $1$ through $3$ are $4$-wise crossings and types 
$4$ through $21$ are $3$-wise crossings. The bundles 
intersecting at each of these $21$ types of vertices are 
listed in Table~\ref{tab:L_j}. For $j=1,2,\ldots,21$, let 
$w_j$ denote the weighted area containing all the crossings 
of type $j$; where the weight is the number of crossings 
per grid cell. To complete the estimates of $\lambda_i(m)$ 
for $i=3,4$, we calculate $w_j$ for all $j$ from the bundles 
intersecting at type $j$ crossings. The values are listed  
in Table~\ref{tab:hex_w_j}. Observe that for two parallel 
strips $\Gamma_i$ and $\Gamma_j$, the area of their intersection 
is $\area(\Gamma_i \cap \Gamma_j) = \area(P(2i-1, 2i, 2j-1, 2j))$;
recall that $P(2i-1, 2i, 2j-1, 2j)$ denotes the parallelogram 
made by the two pairs of parallel lines $\ell_{2i-1}, \ell_{2i}$ 
and $\ell_{2j-1}, \ell_{2j}$, respectively. 

\begin{table}[htpb]
\begin{center}
\begin{minipage}{0.3\linewidth}
\centering
\begin{tabular}{|r|p{3.5cm}|}
\hline
$j$ & Bundles intersecting at type $j$ vertices \\ \hline \hline
$1$ & $\L_6, \L_{12}, \L_3, \L_9$ \\ \hline 
$2$ & $\L_2, \L_8, \L_{11}, \L_5$ \\ \hline 
$3$ & $\L_{10}, \L_4, \L_1, \L_7$ \\ \hline 
$4$ & $\L_2, \L_7, \L_9$ \\ \hline 
$5$ & $\L_6, \L_{11}, \L_1$ \\ \hline 
$6$ & $\L_{10}, \L_3, \L_5$ \\ \hline 
$7$ & $\L_{12}, \L_5, \L_7$ \\ \hline  
\end{tabular}
\end{minipage}
\hspace*{0.03\linewidth}
\begin{minipage}{0.3\linewidth}
\centering
\begin{tabular}{|r|p{3.5cm}|}
\hline
$j$ & Bundles intersecting at type $j$ vertices \\ \hline \hline
$8$ & $\L_4, \L_{11}, \L_9$ \\ \hline 
$9$ & $\L_8, \L_1, \L_3$ \\ \hline
$10$ & $\L_1, \L_5, \L_9$ \\ \hline  
$11$ & $\L_{11}, \L_3, \L_7$ \\ \hline 
$12$ & $\L_{12}, \L_3, \L_9$ \\ \hline 
$13$ & $\L_4, \L_1, \L_7$ \\ \hline 
$14$ & $\L_8, \L_{11}, \L_5$ \\ \hline  
\end{tabular}
\end{minipage}
\hspace*{0.03\linewidth}
\begin{minipage}{0.3\linewidth}
\centering
\begin{tabular}{|r|p{3.5cm}|}
\hline
$j$ & Bundles intersecting at type $j$ vertices \\ \hline \hline
$15$ & $\L_4, \L_8, \L_{12}$ \\ \hline
$16$ & $\L_6, \L_{12}, \L_3$ \\ \hline 
$17$ & $\L_6, \L_{12}, \L_9$ \\ \hline 
$18$ & $\L_2, \L_8, \L_{11}$ \\ \hline 
$19$ & $\L_2, \L_8, \L_5$ \\ \hline
$20$ & $\L_{10}, \L_4, \L_1$ \\ \hline 
$21$ & $\L_{10}, \L_4, \L_7$ \\ \hline   
\end{tabular}
\end{minipage}
\end{center}
\caption{Bundles intersecting at type $j$ vertices for $j=1,2,\ldots,21$.}
\label{tab:L_j}
\end{table}

$\bullet \hspace{2mm}$ For $\lambda_4(m)$, all the $4$-wise 
crossings that are not at grid vertices, are at the centers 
of the grid cells; we have
\[ w_1= \area(\Gamma_6 \cap \Gamma_{12} \cap \Gamma_3 \cap \Gamma_9) =
\area(\Gamma_3 \cap \Gamma_9) = \area(P(5,6,17,18)) = \frac{\sqrt3}{4}.  
\]
Types $2$ and $3$ are $120^{\circ}$ and $240^{\circ}$ rotations 
of type $1$, respectively; therefore by symmetry, $w_1 = w_2 = w_3$.

Consequently, we have
\[ \lambda_4(n) 
= \frac{a_4 + \sum\nolimits_{j=1}^{3} w_j}{\area(\sigma_0)} - O(m) 
= \left( \frac{1}{2} + \frac{3}{8} \right)m^2 - O(m)
= \frac{7 m^2}{8} - O(m). \]

Lastly, we estimate $\lambda_3(m)$. Besides $3$-wise crossings 
at grid vertices in $H$ (whose number is proportional to $a_3$), 
there are $18$ types of $3$-wise crossings \ie, types $4$ through 
$21$, on the boundary or in the interior of the grid cells in $H$. 

\begin{itemize} \itemsep 1pt
\item For types $4,5$, and $6$, there are two crossings per grid 
cell; and 
\begin{align*}
w_4 &= 2 \cdot \area(\Gamma_2 \cap \Gamma_7 \cap \Gamma_9)
= 2 \cdot (\area(P(3,4,17,18)) - \area(1,13,17) - \area(4,14,18)) \\
&= 2 \cdot \left(\frac{2}{\sqrt3} - \frac{1}{4\sqrt3} - \frac{1}{4\sqrt3} \right) = \sqrt3.
\end{align*}
Types $5$ and $6$ are $120^{\circ}$ and $240^{\circ}$ 
rotations of type $4$, respectively; therefore by symmetry, 
$w_4 = w_5 = w_6$.

\item For types $7,8$, and $9$, there are four crossings per grid 
cell; and
\begin{align*}
w_7 &= 4 \cdot \area(\Gamma_{12} \cap \Gamma_5 \cap \Gamma_7)
= 4 \cdot (\area(P(9,10,13,14)) - \area(10,13,23) - \area(9,14,24)) \\
&= 4 \cdot \left(\frac{2\sqrt3}{5} - \frac{\sqrt3}{20} - \frac{\sqrt3}{20} \right)
= \frac{6\sqrt3}{5}. 
\end{align*}
Types $8$ and $9$ are $120^{\circ}$ and $240^{\circ}$ 
rotations of type $7$, respectively; therefore by symmetry, 
$w_7 = w_8 = w_9$.

\item For types $10,11$, there are six crossings per grid cell; and
\begin{align*}
w_{10} &= 6 \cdot \area(\Gamma_1 \cap \Gamma_5 \cap \Gamma_9)
= 6 \cdot (\area(P(1,2,17,18)) - \area(1,9,17) - \area(2,10,18)) \\
&= 6 \cdot \left(\frac{2\sqrt3}{7} - \frac{\sqrt3}{28} - \frac{\sqrt3}{28} \right)
= \frac{9\sqrt3}{7}. 
\end{align*}
Type $11$ is the reflection in a vertical line of type $10$; therefore 
by symmetry, $w_{10} = w_{11}$. 

\item For types $12,13$, and $14$, there are two crossings per grid 
cell; and
\[ w_{12} = 2 \cdot \area(\Gamma_{12} \cap \Gamma_3 \cap \Gamma_9) =
2 \cdot \area(\Gamma_3 \cap \Gamma_9) = 2 \cdot \area(P(5,6,17,18))
= \frac{\sqrt3}{2}. \]
Types $13$ and $14$ are $120^{\circ}$ and $240^{\circ}$ 
rotations of type $12$, respectively; therefore by symmetry,  
$w_{12} = w_{13} = w_{14}$. 

\item For type $15$, there are two crossings per grid cell; and 
\begin{align*}
w_{15} &= 2\cdot \area(\Gamma_4 \cap \Gamma_8 \cap \Gamma_{12})
= 2 \cdot (\area(P(15,16,23,24)) - \area(7,15,24) - \area(8,16,23)) \\
&= 2 \cdot \left(\frac{2}{\sqrt3} - \frac{1}{4\sqrt3} - \frac{1}{4\sqrt3} \right) = \sqrt3.
\end{align*}

\item For types $16$ through $21$, there is one crossing per 
grid cell; and
\begin{align*}
w_{16} &= \area(\Gamma_6 \cap \Gamma_{12} \cap \Gamma_3 - \Gamma_9) 
= \area(\Gamma_{12} \cap \Gamma_3) - \area(\Gamma_{12} \cap \Gamma_3 \cap \Gamma_9) \\
&= \area(P(5,6,23,24)) - \area(P(5,6,17,18))
= \frac{\sqrt3}{2} - \frac{\sqrt3}{4} = \frac{\sqrt3}{4}.
\end{align*}
Type $17$ is the reflection in a vertical line of type $16$, 
types $18$ and $20$ are $120^{\circ}$ and $240^{\circ}$ 
rotations of type $16$, respectively. Types $19$ and $21$ are 
$120^{\circ}$ and $240^{\circ}$ rotations of type $17$, respectively. 
Therefore by symmetry, $w_{16} = w_{17} = w_{18} = w_{19} = w_{20} = w_{21}$.
\end{itemize}

Consequently, we have
\[ \lambda_3(n) 
= \frac{a_3 + \sum\nolimits_{j=4}^{21} w_j}{\area(\sigma_0)} - O(m) 
= \left( \frac32 + \frac32 + \frac95 + \frac97 + \frac34 + \frac12 + \frac34 \right)m^2 - O(m)
= \frac{283}{35}m^2 - O(m). \]
\renewcommand*{\arraystretch}{1.3}
\begin{table}[H]
\begin{center}
\begin{tabular}{|c||c|c|c|c|c|c|c|c|c|c|c|}
\hline 
$j$ & $1$ & $2$ & $3$ & $4$ & $5$ & $6$ & $7$ & $8$ & $9$ & $10$ & $11$ \\ \hline 
$w_j$ & $\frac{\sqrt3}{4}$ & $\frac{\sqrt3}{4}$ & $\frac{\sqrt3}{4}$ 
& $\sqrt3$ & $\sqrt3$ & $\sqrt3$
& $\frac{6\sqrt3}{5}$ & $\frac{6\sqrt3}{5}$ & $\frac{6\sqrt3}{5}$ 
& $\frac{9\sqrt3}{7}$ & $\frac{9\sqrt3}{7}$ \\ \hline 
\multicolumn{12}{c}{ } \\ \cline{1-11} 
$j$ & $12$ & $13$ & $14$ & $15$ & $16$ & $17$ & $18$ & $19$ & $20$ & $21$ \\ \cline{1-11} 
$w_j$ & $\frac{\sqrt3}{2}$ & $\frac{\sqrt3}{2}$ & $\frac{\sqrt3}{2}$ 
& $\sqrt3$ & $\frac{\sqrt3}{4}$ & $\frac{\sqrt3}{4}$
& $\frac{\sqrt3}{4}$ & $\frac{\sqrt3}{4}$ & $\frac{\sqrt3}{4}$ 
& $\frac{\sqrt3}{4}$ \\ \cline{1-11} 
\end{tabular}
\caption{Values of $w_j$ for $j=1, \ldots, 21$.}
\label{tab:hex_w_j}
\end{center}
\end{table}
\renewcommand*{\arraystretch}{1}

The values of $\lambda_i(m)$, for $i=3,\ldots,12$, are summarized 
in Table~\ref{tab:lambda-12hex}; for convenience the linear terms 
are omitted. Since $m = n/12$, $\lambda_i$ can be also viewed as 
a function of $n$.
\renewcommand*{\arraystretch}{1.3}
\begin{table}[H]
\begin{center}
\begin{tabular}{|c||c|c|c|c|c|c|c|c|c|c|}
\hline 
$i$ & $3$ & $4$ & $5$ & $6$ & $7$ & $8$ & $9$ & $10$ & $11$ & $12$ \\ \hline 
$\lambda_i(m)$ & $\frac{283 m^2}{35}$ & $\frac{7 m^2}{8}$ & $\frac{7 m^2}{20}$ 
& $\frac{27 m^2}{140}$ & $\frac{23 m^2}{140}$ & $\frac{19 m^2}{280}$
& $\frac{9 m^2}{140}$ & $\frac{13 m^2}{280}$ & $\frac{m^2}{70}$ 
& $\frac{m^2}{10}$ \\ \hline 
$\lambda_i(n)$ & $\frac{283 n^2}{35 \cdot 144}$ & $\frac{7 n^2}{8 \cdot 144}$ 
& $\frac{7 n^2}{20 \cdot 144}$ & $\frac{27 n^2}{140 \cdot 144}$ 
& $\frac{23 m^2}{280 \cdot 144}$ & $\frac{19 n^2}{280 \cdot 144}$
& $\frac{9 n^2}{140 \cdot 144}$ & $\frac{13 n^2}{280 \cdot 144}$ 
& $\frac{n^2}{70 \cdot 144}$ & $\frac{n^2}{10 \cdot 144}$ \\ \hline
\end{tabular}
\caption{The asymptotic values of $\lambda_i(m)$ and $\lambda_i(n)$ for $i=3,\ldots,12$.}
\label{tab:lambda-12hex}
\end{center}
\end{table}
\renewcommand*{\arraystretch}{1}

The multiplicative factor in Eq.\,\eqref{eq:recurrence-k}
is bounded from below as follows:   
\begin{align*} 
F(n) 
& \geq \prod_{i=3}^{12} B_i^{\lambda_i(n)} 
\geq 2^{\frac{283 n^2}{35 \cdot 144}} 
\cdot 8^{\frac{7 n^2}{8 \cdot 144}} 
\cdot 62^{\frac{7 n^2}{20 \cdot 144}} 
\cdot 908^{\frac{27 n^2}{140 \cdot 144}} 
\cdot 24698^{\frac{23 n^2}{140 \cdot 144}} 
\cdot 1232944^{\frac{19 n^2}{280 \cdot 144}} \\
& \cdot 112018190^{\frac{9 n^2}{140 \cdot 144}} 
\cdot 18410581880^{\frac{13 n^2}{280 \cdot 144}} 
\cdot 5449192389984^{\frac{n^2}{70 \cdot 144}} 
\cdot 2894710651370536^{\frac{n^2}{10 \cdot 144}} 
\cdot 2^{-O(n)}. 
\end{align*}

\begin{figure}[ht]
\centering
\includegraphics[scale=0.7]{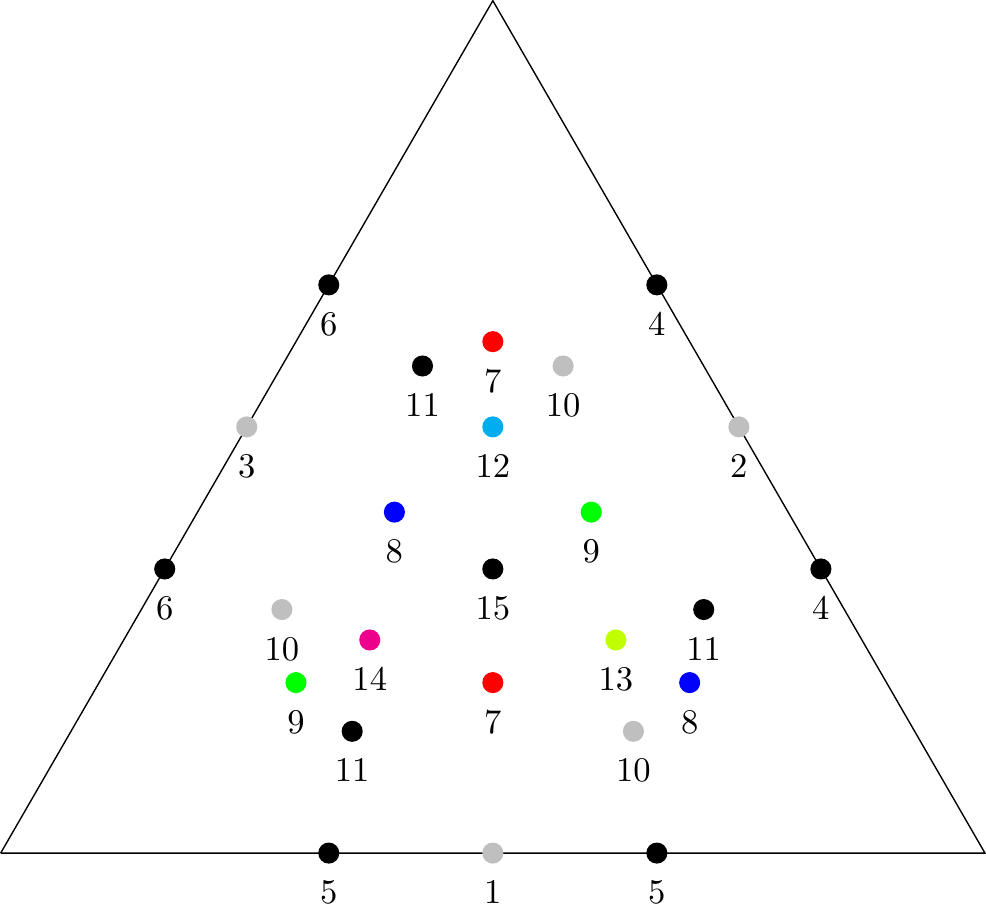}
\caption{In the $12$-gon in the middle of $H$, all the triangular 
grid cells contain $3$-crossings and $4$-crossings of all types $1$ 
through $15$. In other grid cells of the construction only some of 
these types appear.}
\label{fig:f18}
\end{figure}

We prove by induction on $n$ that $T(n) \geq 2^{cn^2 - O(n \log{n})}$ 
for a suitable constant $c>0$. It suffices to choose $c$ (using 
the values of $B_i$ for $i=3,\ldots,12$ in Table~\ref{tab:A_nB_n}) 
so that 
\begin{align*} 
\frac{1}{144} 
& \Big( \frac{283}{35} + \frac{7}{8} \log{8}  + \frac{7}{20} \log{62}
+ \frac{27}{140} \log{908} + \frac{23}{140} \log{24698} \\ 
& + \frac{19}{280} \log{1232944} + \frac{9}{140} \log{112018190} + \frac{13}{280} \log{18410581880} \\ 
& + \frac{1}{70} \log{5449192389984} + \frac{1}{10} \log{2894710651370536} \Big) \geq \frac{11c}{12}. 
\end{align*}

The above inequality holds if we set
\begin{equation} \label{eq:0.2083}
\begin{split}
c & = \frac{1}{132} \Big( \frac{283}{35} + \frac{7}{8} \log{8} + \frac{7}{20} \log{62}
+ \frac{27}{140} \log{908} +\frac{23}{140} \log{24698} \\ 
& + \frac{19}{280} \log{1232944} + \frac{9}{140} \log{112018190} + \frac{13}{280} \log{18410581880} \\ 
& + \frac{1}{70} \log{5449192389984} + \frac{1}{10} \log{2894710651370536} \Big) > 0.2083,
\end{split}
\end{equation}
and the lower bound in Theorem~\ref{thm:main} follows. \qed

\section{Conclusion} \label{sec:conclusion}

We analyzed several recursive constructions derived from arrangements 
of lines with $3$, $4$, $6$, $8$, and $12$ distinct slopes; in 
two different styles (rectangular and hexagonal). The hexagonal 
construction with $12$ slopes yields the lower bound $b_n \geq 0.2083\, n^2$ 
for large $n$. We think that increasing the number of slopes will further
increase the lower bound, and likely the proof complexity at the same time.
The questions of how far can one go and whether there are
other more efficient variants remain. We conclude with the following questions.

\begin{enumerate} \itemsep 1pt
\item What lower bounds on $B_n$ can be deduced from line 
arrangements with a higher number of slopes? In particular, 
hexagonal and rectangular constructions with $16$ slopes seem 
to be the most promising candidates. Note that the value of 
$B_{16}$ is currently unknown. 

\item What lower bounds on $B_n$ can be obtained from rhombic 
tilings of a centrally symmetric octagon? Or from those of a 
centrally symmetric $k$-gon for some other even $k \geq 10$?
No closed formulas for the number of such tilings seem to be 
available at the time of the present writing. However, suitable 
estimates could perhaps be deduced from previous results; see, 
\eg,~\cite{DMB01,DMB04,Eln97,HW15}.
\end{enumerate}

\appendix

\section{Rectangular construction with $8$ slopes} \label{sec:rectangular8}

\begin{figure}[htpb]
\centering
\includegraphics[scale=0.45]{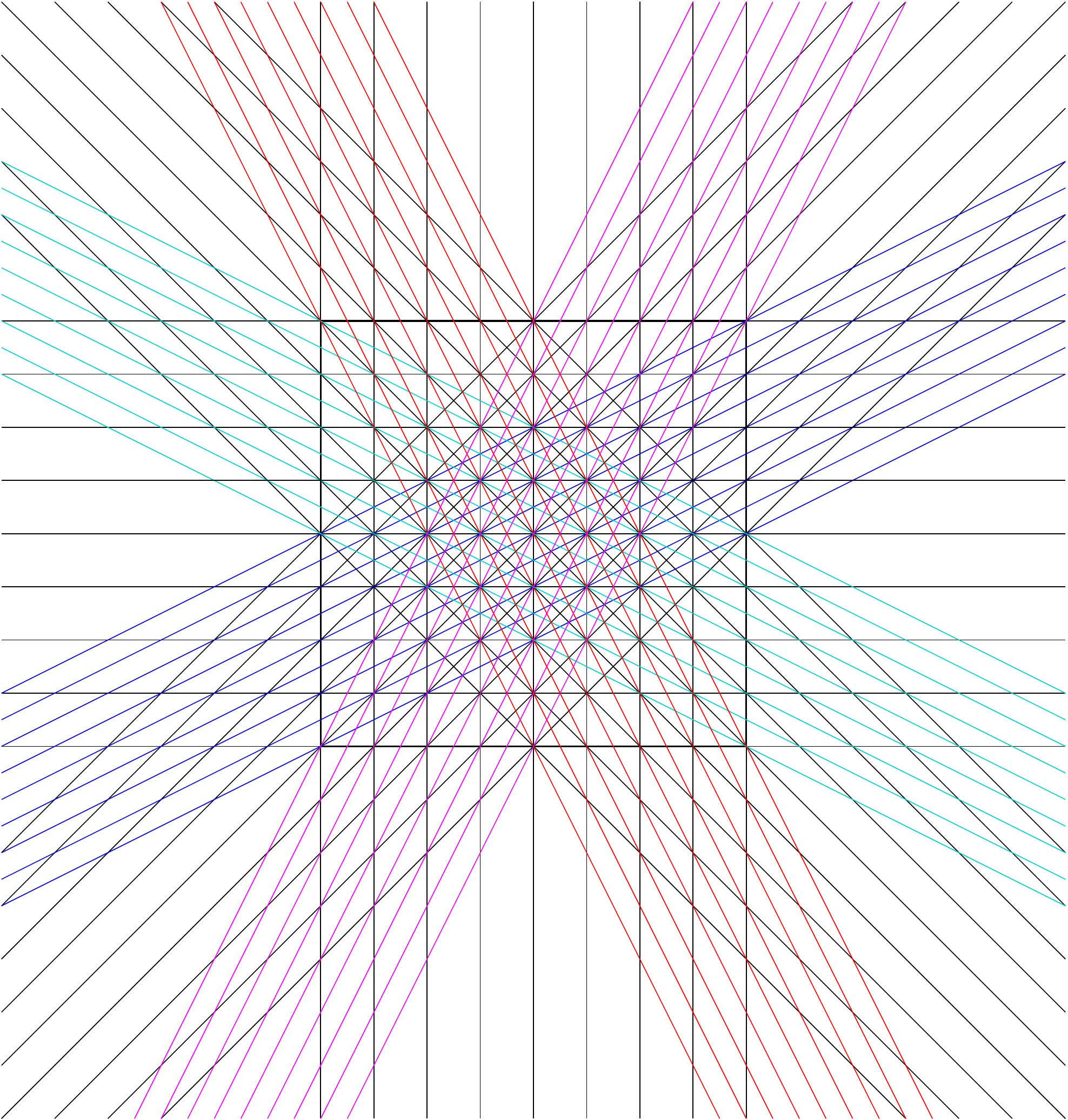}
\caption{Construction with $8$ slopes.}
\label{fig:f6}
\end{figure}

We describe and analyze a rectangular construction with lines 
of $8$ slopes. See Fig.~\ref{fig:f6}. Consider $8$ bundles of 
parallel lines whose slopes are $0, \infty, \pm1/2, \pm1, \pm2$. 
The axes of all parallel strips are all incident to the center 
of $U$. This construction yields the lower bound 
$b_n \geq 0.1999\, n^2$ for large $n$.

\begin{figure}[t]
\hspace{0.05\linewidth}
\begin{minipage}[b]{0.6\linewidth}
\centering
\includegraphics[width=\linewidth]{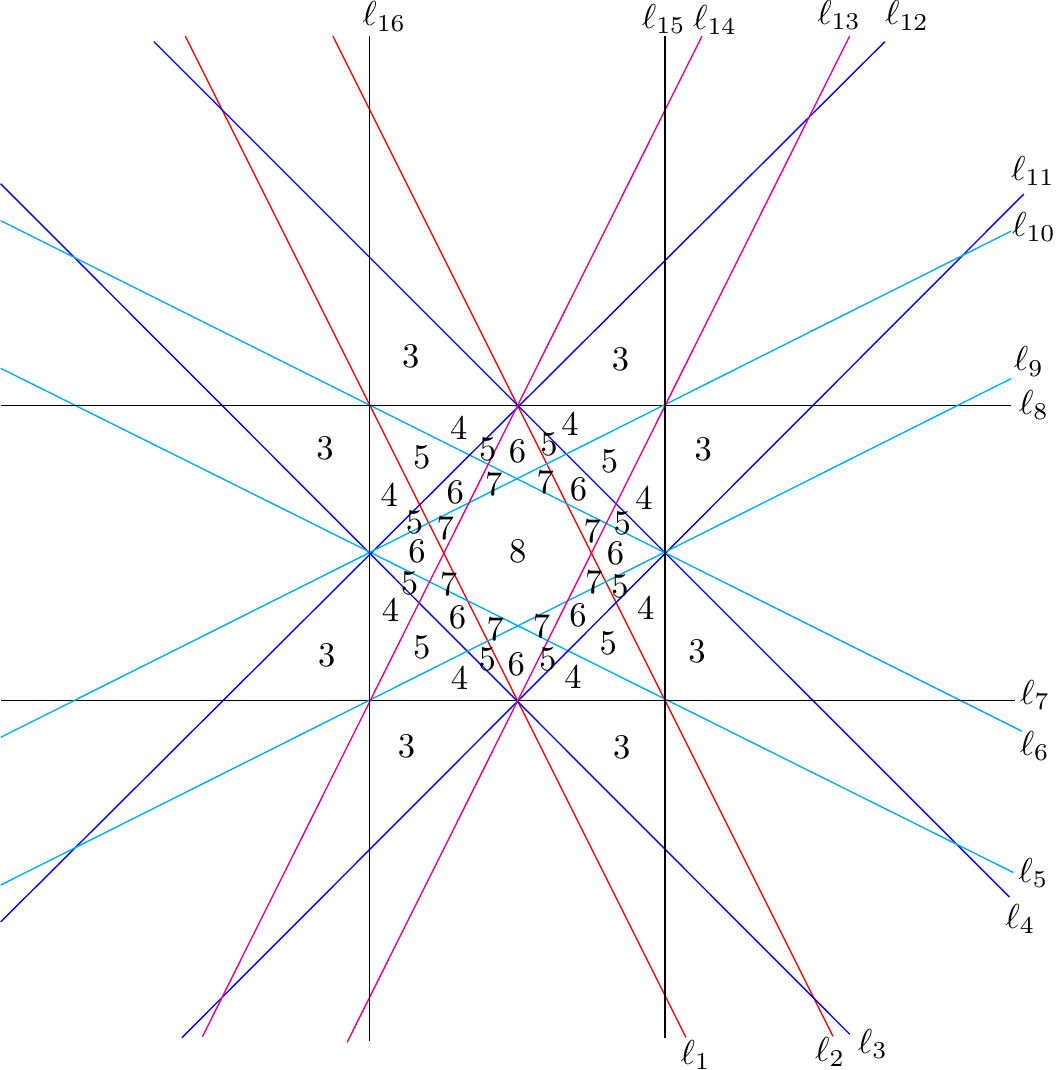}
\vspace{-4.7cm}
\end{minipage}
\hspace{0.05\linewidth}
\begin{minipage}[b]{0.2\linewidth}
\centering
\begin{tabular}{|r||c|c|c|} \hline
$i$ & $\alpha_i$ & $\beta_i$ & $\gamma_i$ \\ \hline \hline
$1$ & $2$ & $1$ & $-1$ \\ \hline
$2$ & $2$ & $1$ & $-2$ \\ \hline
$3$ & $1$ & $1$ & $-0.5$ \\ \hline 
$4$ & $1$ & $1$ & $-1.5$ \\ \hline 
$5$ & $1$ & $2$ & $-1$ \\ \hline 
$6$ & $1$ & $2$ & $-2$ \\ \hline 
$7$ & $0$ & $1$ & $0$ \\ \hline 
$8$ & $0$ & $1$ & $-1$ \\ \hline 
$9$ & $1$ & $-2$ & $0$ \\ \hline 
$10$ & $1$ & $-2$ & $1$ \\ \hline
$11$ & $1$ & $-1$ & $-0.5$ \\ \hline 
$12$ & $1$ & $-1$ & $0.5$ \\ \hline  
$13$ & $2$ & $-1$ & $-1$ \\ \hline 
$14$ & $2$ & $-1$ & $0$ \\ \hline 
$15$ & $1$ & $0$ & $-1$ \\ \hline 
$16$ & $1$ & $0$ & $0$ \\ \hline 
\end{tabular}
\end{minipage}
\caption{Left: The eight parallel strips and the corresponding 
covering multiplicities. These numbers only reflect incidences 
at the grid vertices made by axis-aligned lines. 
Right: Coefficients of the lines $\ell_i$ for $i=1,2,\ldots,16$.}
\label{fig:f7}
\end{figure}

Let $\L=\L_1 \cup \ldots \cup \L_8$ be the partition of $\L$ 
into eight bundles. The $m$ lines in $\L_i$ are contained in 
the parallel strip $\Gamma_i$ bounded by the two lines $\ell_{2i-1}$ 
and $\ell_{2i}$, for $i=1,\ldots,8$. The equation of line 
$\ell_i$ is $\alpha_i x +\beta_i y +\gamma_i=0$, with 
$\alpha_i,\beta_i,\gamma_i$, for $i=1,\ldots,16$ given in 
Fig~\ref{fig:f7}\,(right). Observe that $U = \Gamma_4 \cap \Gamma_8$. 

We refer to lines in $\L_4 \cup \L_8$ (\ie, axis-aligned lines) 
as the \emph{primary} lines, and to rest of the lines as 
\emph{secondary} lines. We refer to the intersection points of 
the primary lines as \emph{grid vertices}. The slopes of 
the primary lines are in $\{0,\infty\}$. The slopes of the 
secondary lines are in $\{\pm1/2,\pm1,\pm2\}$. Note that
\begin{itemize} \itemsep 1pt
\item the distance between consecutive lines in $\L_4$ or $\L_8$ is 
$\frac{1}{m} \left(1-O\left(\frac{1}{m}\right)\right)$;

\item the distance between consecutive lines in $\L_2$ or $\L_6$ is 
$\frac{1}{m \sqrt2} \left(1-O\left(\frac{1}{m}\right)\right)$;

\item the distance between consecutive lines in $\L_1$, $\L_3$, 
$\L_5$ or $\L_7$ is $\frac{1}{m \sqrt5} \left(1-O\left(\frac{1}{m}\right)\right)$.
\end{itemize}

Let $\sigma_0 = \sigma_0(m)$ denote the basic parallelogram (here, 
square) determined by consecutive axis-aligned lines (\ie, lines 
in $\L_4 \cup \L_8$); the side length of $\sigma_0$ is 
$\frac{1}{m} \left(1-O\left(\frac{1}{m}\right)\right)$. 
We refer to these basic parallelograms as \emph{grid cells}. 
Let $U'$ be the smaller square made by 
$\ell_5,\ell_6,\ell_{13},\ell_{14}$, 
\ie, $U' = \Gamma_3 \cap \Gamma_7$; the similarity ratio 
$\rho(U',U)$ is equal to $\frac{1}{\sqrt5}$. We have   
\begin{align*}
\area(U) &= 1, \\
\area(U') &= \frac{\area(U)}{5}=\frac{1}{5}, \\
\area(\sigma_0) &= \frac{1}{m^2}\left(1-O\left(\frac{1}{m}\right)\right). 
\end{align*}

For $i=3,4,\ldots,8$, let $a_i$ denote the area of the (not 
necessarily connected) region covered by exactly $i$ of the 
$8$ strips. Recall that $\area(i,j,k)$ denotes the area of 
the triangle made by $\ell_i$, $\ell_j$ and $\ell_k$. We have
\begin{align*}
a_3 &= 8 \cdot \area(3,7,15) = 1, \\ 
a_4 &= 8 \cdot \area(5,7,11) = \frac13, \\ 
a_5 &= 4 \left( 2 \cdot \area(5,11,13) + \area(2,5,11) \right) = \frac{7}{30}, \\ 
a_6 &= 4 \left( \area(6,11,13) - 2 \cdot \area(2,11,9) - \area(2,9,13) \right) = \frac15, \\ 
a_7 &= 8 \cdot \area(5,9,13 )= \frac{1}{15}, \\ 
a_8 &= \area(U') - 4 \cdot \area(5,9,13) = \frac15 - \frac{1}{30} = \frac16. 
\end{align*}

Observe that $a_4+a_5+a_6+a_7+a_8 = \area(U) = 1$. Recall 
that $\lambda_i(m)$ denote the number of $i$-wise crossings 
where each bundle consists of $m$ lines. Note that 
$\lambda_i(m)$ is proportional to $a_i$, for $i=4,5,6,7,8$. 
Indeed, $\lambda_i(m)$ is equal to the number of grid points 
that lie in a region covered by $i$ parallel strips, which 
is roughly equal to the ratio $\frac{a_i}{\area(\sigma_0)}$, 
for $i=4,5,6,7,8$. More precisely, taking also the boundary 
effect of the relevant regions into account, we obtain
\begin{align*}
\lambda_4(m) &= \frac{a_4}{\area(\sigma_0)} - O(m) = \frac{m^2}{3} - O(m), \\
\lambda_5(m) &= \frac{a_5}{\area(\sigma_0)} - O(m) = \frac{7 m^2}{30} - O(m), \\
\lambda_6(m) &= \frac{a_6}{\area(\sigma_0)} - O(m) = \frac{m^2}{5} - O(m), \\
\lambda_7(m) &= \frac{a_7}{\area(\sigma_0)} - O(m) = \frac{m^2}{15} - O(m), \\
\lambda_8(m) &= \frac{a_8}{\area(\sigma_0)} - O(m) = \frac{m^2}{6} - O(m). 
\end{align*}

\begin{figure}[htpb]
\begin{minipage}[b]{0.6\linewidth}
\centering
\includegraphics[scale=1]{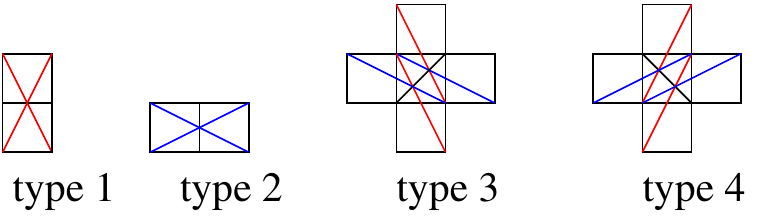}
\vspace{-1.5cm}
\end{minipage}
\hspace{0.01\linewidth}
\begin{minipage}[b]{0.35\linewidth}
\centering
\begin{tabular}{|r|p{3.5cm}|} \hline
$j$ & Bundles intersecting at vertices of type $j$ \\ \hline \hline
$1$ & $\L_4, \L_1, \L_7$ \\ \hline 
$2$ & $\L_8, \L_3, \L_5$ \\ \hline 
$3$ & $\L_6, \L_1, \L_3$ \\ \hline 
$4$ & $\L_2, \L_5, \L_7$ \\ \hline 
\end{tabular}
\end{minipage}
\caption{Left: Other types of $3$-wise crossings. 
Right: Intersecting bundles for these crossings.}
\label{fig:f3}
\end{figure}

For estimating $\lambda_3(m)$, in addition to considering 
$3$-wise crossings in the exterior of $U$, we also observe 
$3$-wise crossings on the boundaries or in the interior of 
the small grid cells contained in some regions of $U$. 
Specifically, we distinguish exactly four types of $3$-wise 
crossings, as illustrated and specified in Fig.~\ref{fig:f3}. 
For $j=1,2,3,4$, let $w_j$ denote the weighted area 
containing all crossings of type $j$, where the weight is 
the number of $3$-wise crossings per grid cell. To complete 
the estimate of $\lambda_3(m)$, we calculate $w_j$ for all $j$, 
from the bundles intersecting at crossings of type $j$; listed 
in Fig.~\ref{fig:f3}\,(right).

Recall that for two parallel strips $\Gamma_i$ and $\Gamma_j$, 
the area of their intersection is 
$\area(\Gamma_i \cap \Gamma_j) = \area(P(2i-1, 2i, 2j-1, 2j))$;
where $P(2i-1, 2i, 2j-1, 2j)$ denotes the parallelogram bounded by
the two pairs of parallel lines $\ell_{2i-1}, \ell_{2i}$ and 
$\ell_{2j-1}, \ell_{2j}$, respectively. For types $1$ and $2$, 
there is one crossing per grid cell and for types $3$ and $4$, 
there are two crossings per grid cell. Therefore we have,
\begin{itemize}\itemsep 1pt
\item $w_1 = \area(\Gamma_4 \cap \Gamma_1 \cap \Gamma_7)
  = \area(\Gamma_1 \cap \Gamma_7) = \area(P(1,2,13,14)) = 1/4$, 
\item $w_2 = \area(\Gamma_8 \cap \Gamma_3 \cap \Gamma_5)
  = \area(\Gamma_3 \cap \Gamma_5) = \area(P(5,6,9,10)) = 1/4$, 
\item $w_3 = 2 \cdot \area(\Gamma_1 \cap \Gamma_3 \cap \Gamma_6)
  = 2 \cdot (\area(P(1,2,5,6)) - 2 \cdot \area(2,5,11)) 
= 2 \cdot (1/3 - 1/12) = 1/2$, 
\item $w_4 = 2 \cdot \area(\Gamma_2 \cap \Gamma_5 \cap \Gamma_7) = 1/2$. 
\end{itemize} 

It follows that
\[ \lambda_3(m) = \frac{a_3 + \sum_{j=1}^4 w_j}{\area(\sigma_0)} - O(m)=
\left(1 +\frac14 +\frac14 +\frac12 +\frac12 \right) m^2 -O(m) = \frac{5 m^2}{2} - O(m). \]

The values of $\lambda_i(m)$, for $i=3,4,\ldots,8$, are summarized 
in Table~\ref{tab:lambda-8}; for convenience the linear terms 
are omitted. Since $m = n/8$, $\lambda_i$ can be also viewed 
as a function of $n$. 
\renewcommand*{\arraystretch}{1.2}
\begin{table}[ht]
\begin{center}
\begin{tabular}{|c||c|c|c|c|c|c|}
\hline
$i$ & $3$ & $4$ & $5$ & $6$ & $7$ & $8$ \\
\hline 
$\lambda_i(m)$ & $\frac{5 m^2}{2}$ & $\frac{m^2}{3}$ & $\frac{7 m^2}{30}$ 
& $\frac{m^2}{5}$ & $\frac{m^2}{15}$ & $\frac{m^2}{6}$ \\ 
\hline 
$\lambda_i(n)$ & $\frac{5 n^2}{2 \cdot 64}$ & $\frac{n^2}{3 \cdot 64}$
& $\frac{7 n^2}{30 \cdot 64}$ & $\frac{n^2}{5 \cdot 64}$ 
& $\frac{n^2}{15 \cdot 64}$ & $\frac{n^2}{6 \cdot 64}$ \\
\hline
\end{tabular}
\caption{The asymptotic values of $\lambda_i(m)$ and $\lambda_i(n)$ 
for $i=3,4,\ldots,8$.}
\label{tab:lambda-8}
\end{center}
\end{table}
\renewcommand*{\arraystretch}{1}

The multiplicative factor in Eq.\,\eqref{eq:recurrence-k}
is bounded from below as follows:  
\[ F(n) \geq \prod_{i=3}^8 B_i^{\lambda_i(n)} \geq
2^{\frac{5 n^2}{2 \cdot 64}} \cdot 8^{\frac{n^2}{3 \cdot 64}} \cdot 
62^{\frac{7 n^2}{30 \cdot 64}} \cdot 908^{\frac{n^2}{5 \cdot 64}} \cdot 
24698^{\frac{n^2}{15 \cdot 64}} \cdot 1232944^{\frac{n^2}{6 \cdot 64}} \cdot 2^{-O(n)}. \]

We prove by induction on $n$ that $T(n) \geq 2^{cn^2 -O(n \log{n})}$ 
for a suitable constant $c>0$. It suffices to choose $c$ (using the 
values of $B_i$ for $i=3,\ldots,8$ in Table~\ref{tab:A_nB_n}) so that
\[ \frac{1}{64} \left( \frac{5}{2} + \frac13 \log{8} + \frac{7}{30} \log{62} 
+ \frac15 \log{908} + \frac{1}{15} \log{24698} + \frac16 \log{1232944} \right)
\geq \frac{7c}{8}. \]

The above inequality holds if we set
\begin{equation} \label{eq:0.1999}
c = \frac{1}{56} \left(\frac{5}{2} + 1 + \frac{7}{30} \log{62} 
+ \frac15 \log{908} + \frac{1}{15} \log{24698} 
+ \frac16 \log{1232944} \right) > 0.1999,
\end{equation}
and this yields the lower bound $B_n \geq 2^{cn^2 - O(n \log{n})}$, 
for some constant $c > 0.1999$. In particular, we have 
$B_n \geq 2^{0.1999\, n^2}$ for large $n$.

\section{Rectangular construction with $12$ slopes} \label{sec:rectangular12}

\begin{figure}[H]
\centering
\includegraphics[scale=0.8]{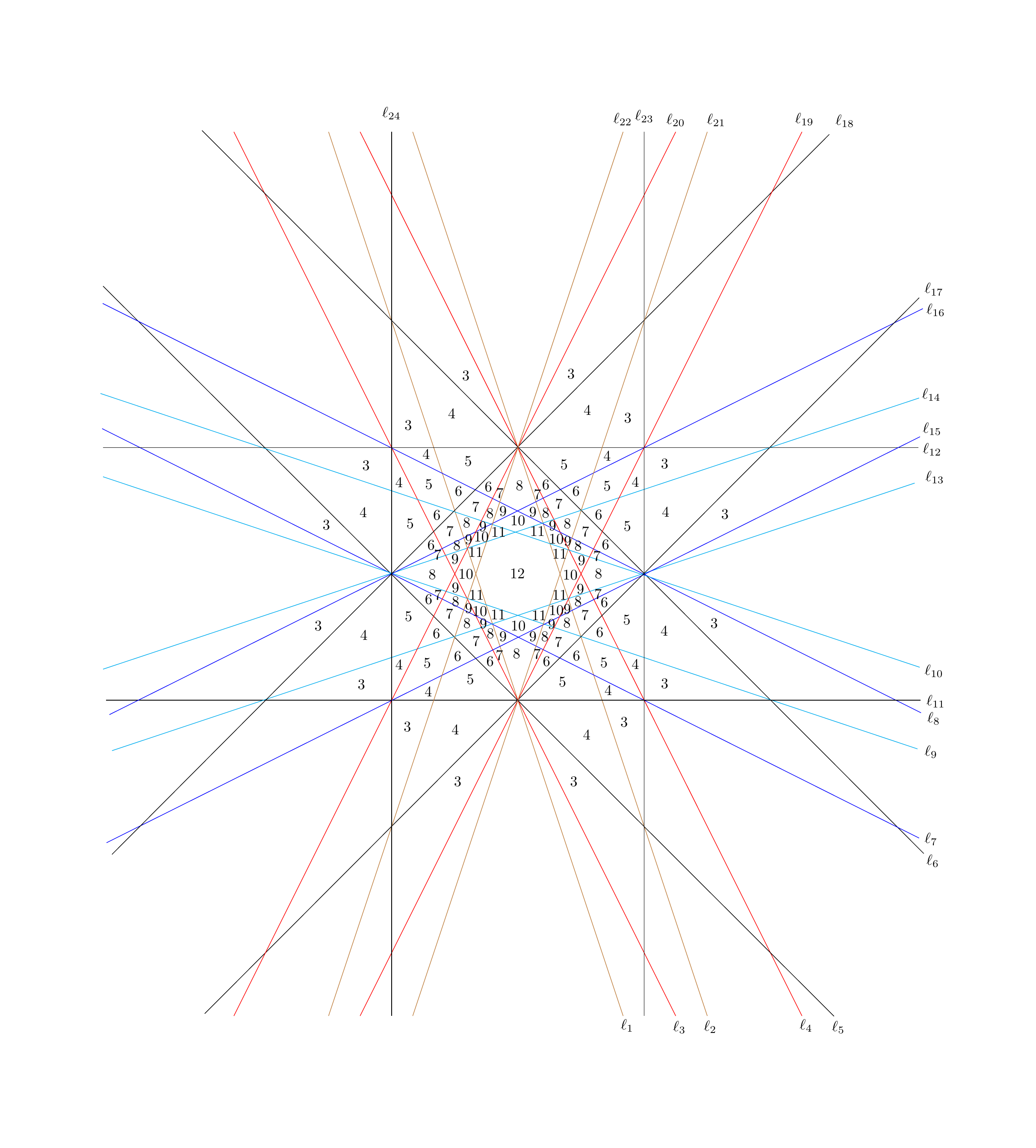}
\caption{Construction with $12$ slopes. The twelve parallel 
strips and the corresponding covering multiplicities. These 
numbers only reflect incidences at the grid vertices made 
by axis-aligned lines.}
\label{fig:f14}
\end{figure}

We next describe and analyze a rectangular construction with 
lines of $12$ slopes. Consider $12$ bundles of parallel lines 
whose slopes are $0,\infty,\pm1/3,\pm1/2,\pm1,\pm2,\pm3$. 
The axes of all parallel strips are all incident to the center 
of $U =[0,1]^2$. Refer to Fig.~\ref{fig:f14}. This construction 
yields the lower bound $b_n \geq 0.2053\, n^2$ for large $n$. 

\begin{table}[htpb]
\begin{center}
\begin{minipage}{0.32\linewidth}
\centering
\begin{tabular}{|r||c|c|c|} \hline
$i$ & $\alpha_i$ & $\beta_i$ & $\gamma_i$ \\ \hline \hline
$1$ & $3$ & $1$ & $-1.5$ \\ \hline
$2$ & $3$ & $1$ & $-2.5$ \\ \hline
$3$ & $2$ & $1$ & $-1$ \\ \hline 
$4$ & $2$ & $1$ & $-2$ \\ \hline 
$5$ & $1$ & $1$ & $-0.5$ \\ \hline 
$6$ & $1$ & $1$ & $-1.5$ \\ \hline
$7$ & $1$ & $2$ & $-1$ \\ \hline 
$8$ & $1$ & $2$ & $-2$ \\ \hline
\end{tabular}
\end{minipage}
\begin{minipage}{0.32\linewidth}
\centering
\begin{tabular}{|r||c|c|c|} \hline
$i$ & $\alpha_i$ & $\beta_i$ & $\gamma_i$ \\ \hline \hline
$9$ & $1$ & $3$ & $-1.5$ \\ \hline 
$10$ & $1$ & $3$ & $-2.5$ \\ \hline
$11$ & $0$ & $1$ & $0$ \\ \hline 
$12$ & $0$ & $1$ & $-1$ \\ \hline 
$13$ & $-1$ & $3$ & $-0.5$ \\ \hline 
$14$ & $-1$ & $3$ & $-1.5$ \\ \hline 
$15$ & $-1$ & $2$ & $0$ \\ \hline 
$16$ & $-1$ & $2$ & $-1$ \\ \hline
\end{tabular}
\end{minipage}
\begin{minipage}{0.32\linewidth}
\centering
\begin{tabular}{|r||c|c|c|} \hline
$i$ & $\alpha_i$ & $\beta_i$ & $\gamma_i$ \\ \hline \hline 
$17$ & $-1$ & $1$ & $0.5$ \\ \hline 
$18$ & $-1$ & $1$ & $-0.5$ \\ \hline 
$19$ & $-2$ & $1$ & $1$ \\ \hline 
$20$ & $-2$ & $1$ & $0$ \\ \hline
$21$ & $-3$ & $1$ & $1.5$ \\ \hline 
$22$ & $-3$ & $1$ & $0.5$ \\ \hline
$23$ & $-1$ & $0$ & $1$ \\ \hline 
$24$ & $-1$ & $0$ & $0$ \\ \hline
\end{tabular}
\end{minipage}
\end{center}
\caption{Coefficients of the $24$ lines.}
\label{tab:12lines}
\end{table}

Let $\L=\L_1 \cup \ldots \cup \L_{12}$ be the partition of 
$\L$ into twelve bundles. The $m$ lines in $\L_i$ are contained 
in the parallel strip $\Gamma_i$ bounded by the two lines 
$\ell_{2i-1}$ and $\ell_{2i}$, for $i=1,\ldots,12$. The 
equation of line $\ell_i$ is $\alpha_i x +\beta_i y +\gamma_i=0$, 
with $\alpha_i,\beta_i,\gamma_i$, for $i=1,\ldots,24$ given in 
Table~\ref{tab:12lines}. Observe that $U = \Gamma_6 \cap \Gamma_{12}$.

We refer to lines in $\L_6 \cup \L_{12}$ (\ie, axis-aligned 
lines) as the \emph{primary} lines, and to rest of the lines 
as the \emph{secondary} lines. We refer to the intersection 
points of the primary lines as \emph{grid vertices}. The slopes 
of the primary lines are in $\{0,\infty\}$, and the slopes of 
the secondary lines are in $\{\pm1/3,\pm1/2,\pm1,\pm2,\pm3\}$. 

\newpage
Note that
\begin{itemize} \itemsep 1pt
\item the distance between consecutive lines in $\L_6$ or $\L_{12}$ is 
$\frac{1}{m} \left(1-O\left(\frac{1}{m}\right)\right)$;

\item the distance between consecutive lines in $\L_3$ or $\L_9$ is 
$\frac{1}{m \sqrt2} \left(1-O\left(\frac{1}{m}\right)\right)$;

\item the distance between consecutive lines in $\L_2$, $\L_4$, $\L_8$, or $\L_{10}$ is 
$\frac{1}{m \sqrt5} \left(1-O\left(\frac{1}{m}\right)\right)$;

\item the distance between consecutive lines in $\L_1$, $\L_5$, $\L_7$, or $\L_{11}$ is 
$\frac{1}{m \sqrt{10}} \left(1-O\left(\frac{1}{m}\right)\right)$.
\end{itemize}

Let $\sigma_0=\sigma_0(m)$ denote the basic parallelogram (here, 
square) determined by consecutive axis-aligned lines (\ie, lines 
in $\L_6 \cup \L_{12}$); the side length of $\sigma_0$ is 
$\frac{1}{m} \left(1-O\left(\frac{1}{m}\right)\right)$. 
We refer to these basic parallelograms as \emph{grid cells}. 
Let $U_1 = \Gamma_1 \cap \Gamma_7$, be the square made by 
$\ell_1,\ell_2,\ell_{13},\ell_{14}$, and let 
$U_2 = \Gamma_2 \cap \Gamma_8$, be the smaller square made by 
$\ell_3,\ell_4,\ell_{15},\ell_{16}$.
Note that $\rho(U_1,U)=\frac{1}{\sqrt{10}}$ and 
$\rho(U_2,U)=\frac{1}{\sqrt{5}}$. We also have 
\begin{align*}
\area(U) &= 1, \\ 
\area(U_1) &= \frac{\area(U)}{10}=\frac{1}{10}, \\ 
\area(U_2) &= \frac{\area(U)}{5}=\frac{1}{5}, \\ 
\area(\sigma_0) &= \frac{1}{m^2}\left(1-O\left(\frac{1}{m}\right)\right). 
\end{align*}

For $i=3,\ldots,12$, let $a_i$ denote the area of the (not 
necessarily connected) region covered by exactly $i$ of the 
$12$ strips. Recall that $\area(i,j,k)$ denotes the area of 
the triangle bounded by $\ell_i$, $\ell_j$ and $\ell_k$. We have
\begin{align*}
a_3 &= 8 \cdot (\area(2,11,13) + \area(3,5,23))= 8 \cdot \left( \frac18 + \frac{1}{24}\right) = \frac43, \\
a_4 &= 8 \cdot (\area(2,5,11) + \area(2,7,11)) = 8 \cdot \left( \frac{1}{12} + \frac{1}{120} \right)
= \frac{11}{15}, \\
a_5 &= 4 \cdot \left(\area(11,17,23) - 2 \cdot \area(2,7,11) - 2 \cdot \area(2,7,17) \right) 
= 4 \left( \frac18 - \frac{2}{120} - \frac{2}{120} \right) = \frac{11}{30}, \\
a_6 &= 4 \cdot \left( 2 \cdot \area(7,17,19) + 2 \cdot \area(2,7,17) \right) 
= 4 \cdot \left( \frac{2}{120}+ \frac{2}{120}\right) = \frac{2}{15}, \\
a_7 &= 4 \cdot \left( 2 \cdot \area(7,19,21) + 2 \cdot (\area(9,17,19) - \area(7,17,19)) \right) \\ 
&= 4 \cdot \left( \frac{1}{140} + 2 \cdot \left( \frac{1}{56} - \frac{1}{120} \right) \right)
= \frac{11}{105}, \\
a_8 &= 8 \cdot \left( \area(9,19,21) - \area(7,19,21) \right)
+ 4 \cdot \left( \area(2,9,15) - \area(9,15,19) \right) \\ 
&+ 8 \cdot \area(7,21,25) = \frac{13}{105}, \\
a_9 &= 8 \cdot (\area(7,15,21) + \area(9,15,19)) = 8 \cdot \left( \frac{1}{280} + \frac{1}{840} \right) 
= \frac{4}{105}, \\
a_{10} &= 4 \cdot \left( (\area(7,13,15) -\area(9,13,15)) + (\area(13,19,21) - \area(15,19,21)) \right) \\
&= 4 \cdot \left( \left( \frac{1}{40} - \frac{1}{60} \right)
+ \left( \frac{1}{80} - \frac{1}{120} \right) \right)
= 4 \cdot \left( \frac{1}{120} + \frac{1}{240} \right) = \frac{1}{20}, \\
a_{11} &= 8 \cdot \area(2,13,21) =  \frac{8}{240} = \frac{1}{30}, \\
a_{12} &= \area(U_1) - 4 \cdot \area(9,13,21) = \frac{1}{10} - \frac{4}{240} = \frac{1}{12}. 
\end{align*}

Observe that the region whose area is $\sum_{i=4}^{12} a_i$ 
consists of $U$ and $8$ triangles outside $U$. Therefore, 
\[ \sum\nolimits_{i=4}^{12} a_i = \area(U) + 8 \cdot \area(2,5,11) = 1 + 2/3 = 5/3. \]

Recall that $\lambda_i(m)$ denote the number of $i$-wise 
crossings where each bundle consists of $m$ lines. Note that 
$\lambda_i(m)$ is proportional to $a_i$, for $i=7,8,\ldots,12$. 
Indeed, $\lambda_i(m)$ is equal to the number of grid vertices, 
\ie, intersection points of the axis-parallel lines that lie 
in a region covered by $i$ parallel strips, which is roughly 
equal to the ratio $\frac{a_i}{\area(\sigma_0)}$, for 
$i=7,8,\ldots,12$. More precisely, taking also the boundary 
effect of the relevant regions into account, we obtain
\begin{align*}
\lambda_7(m) &= \frac{a_7}{\area(\sigma_0)} -O(m) =\frac{11 m^2}{105} -O(m), \\
\lambda_8(m) &= \frac{a_8}{\area(\sigma_0)} -O(m) =\frac{13 m^2}{105} -O(m), \\
\lambda_9(m) &= \frac{a_9}{\area(\sigma_0)} -O(m) =\frac{4 m^2}{105} -O(m), \\
\lambda_{10}(m) &= \frac{a_{10}}{\area(\sigma_0)} -O(m) =\frac{m^2}{20} -O(m), \\
\lambda_{11}(m) &= \frac{a_{11}}{\area(\sigma_0)} -O(m) =\frac{m^2}{30} -O(m), \\
\lambda_{12}(m) &= \frac{a_{12}}{\area(\sigma_0)} -O(m) =\frac{m^2}{12} -O(m). 
\end{align*}

\begin{figure}[htpb]
\centering
\includegraphics[scale=0.75]{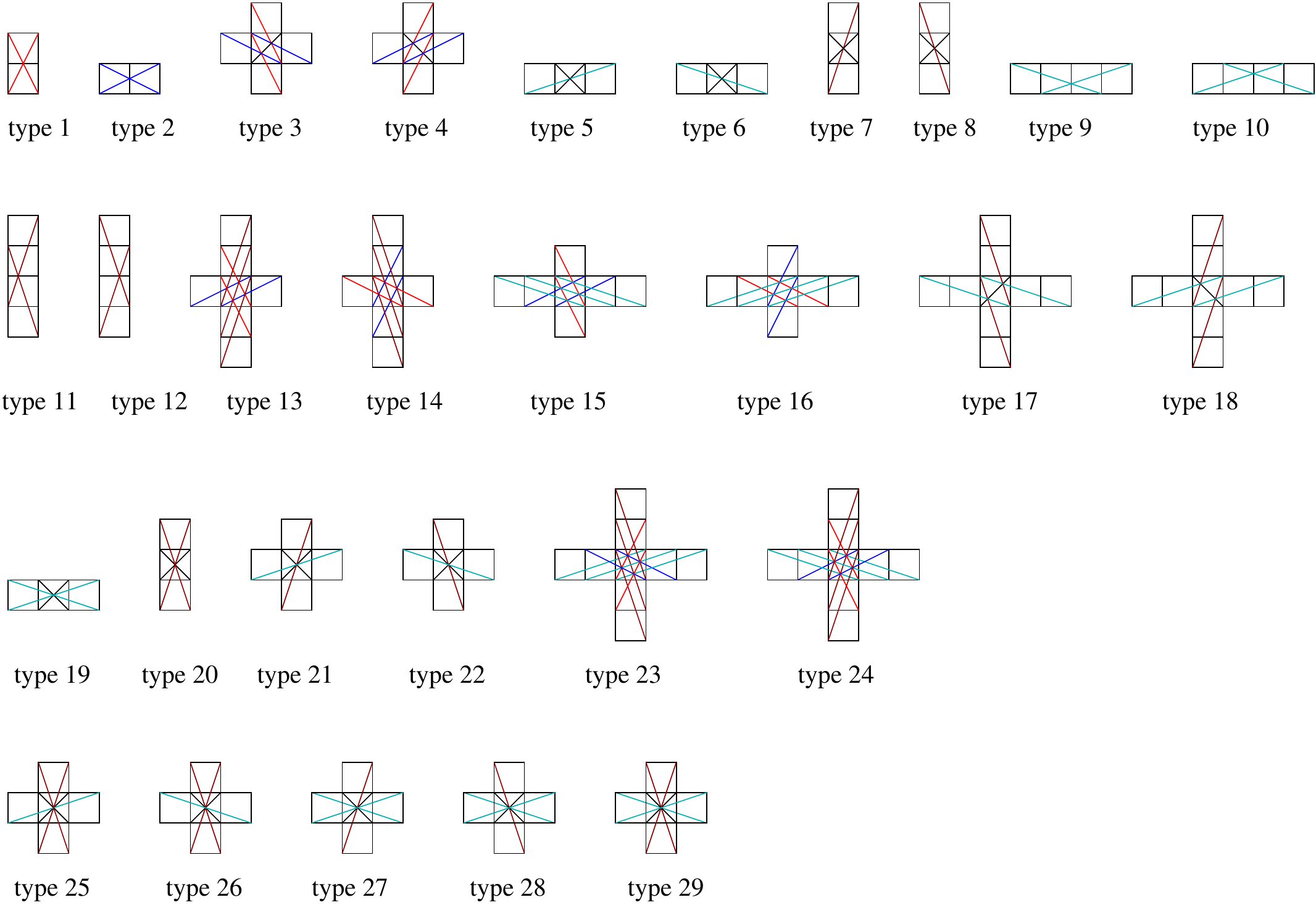}
\caption{Types of incidences of $3$, $4$, $5$, and $6$ lines 
that are not at grid vertices:
$3$-wise crossings: types $1$ through $18$;
$4$-wise crossings: types $19$ through $24$;
$5$-wise crossings: types $25$ through $28$;
$6$-wise crossings: type $29$. \\
For some types, the crossings are in the middle of a grid cell.
To list the coordinates of crossing points, we rescale
the grid cells to the unit square $[0,1]^2$. \\
For types $1$ and $2$, the crossings are at the midpoint of the 
horizontal and the vertical grid edges respectively. 
For type $3$, the crossings are at $(1/3, 1/3)$ and $(2/3, 2/3)$. \\
For type $4$, the crossings are at $(1/3, 2/3)$ and $(2/3, 1/3)$. \\
For types $9$ and $10$, the crossings are on vertical grid edges at height $1/3$ and $2/3$
from the horizontal line below, respectively. \\
For types $11$ and $12$, the crossings are on horizontal grid edges at distance $1/3$ and $2/3$
from the vertical line on the left, respectively. \\
For type $13$, the crossings are at $(1/5, 3/5)$ and $(3/5, 4/5)$ and $(4/5, 2/5)$ and $(2/5, 1/5)$. \\
For type $14$, the crossings are at $(1/5, 2/5)$ and $(2/5, 4/5)$ and $(4/5, 3/5)$ and $(3/5, 1/5)$. \\
For type $15$, the crossings are at $(1/5, 3/5)$ and $(3/5, 4/5)$ and $(4/5, 2/5)$ and $(2/5, 1/5)$. \\
For type $16$, the crossings are at $(1/5, 2/5)$ and $(2/5, 4/5)$ and $(4/5, 3/5)$ and $(3/5, 1/5)$. \\
For type $17$, the crossings are at $(1/4, 1/4)$ and $(3/4, 3/4)$. \\
For type $18$, the crossings are at $(1/4, 3/4)$ and $(3/4, 1/4)$. \\
For type $23$, the crossings are at $(1/5, 2/5)$ and $(2/5, 4/5)$ and $(4/5, 3/5)$ and $(3/5, 1/5)$. \\
For type $24$, the crossings are at $(1/5, 3/5)$ and $(3/5, 4/5)$ and $(4/5, 2/5)$ and $(2/5, 1/5)$. \\
For the other types, the crossings are at $(1/2,1/2)$. 
}
\label{fig:f13}
\end{figure}

For $i=3,4,5,6$, not all the $i$-wise crossings are at 
grid vertices. It can be exhaustively verified (by hand)
that there are $29$ types of such crossings in 
total; see Fig.~\ref{fig:f13}. The bundles intersecting 
at each of these $29$ types of vertices are listed in 
Table~\ref{tab:rec_L_j}. For $j=1,2,\ldots,29$, let 
$w_j$ denote the weighted area containing all crossings of type $j$;
where the weight is the number of crossings per grid cell.
To complete the estimates of $\lambda_i(m)$ for 
$i=3,4,5,6$, we calculate $w_j$ for all $j$ from the bundles 
intersecting at type $j$ crossings. The values are listed 
in Table~\ref{tab:rec_w_j}.

\begin{table}[htpb]
\begin{center}
\begin{minipage}{0.3\linewidth}
\centering
\begin{tabular}{|r|p{3.5cm}|}
\hline
$j$ & Bundles intersecting at type $j$ vertices \\
\hline \hline
$1$ & $\L_2, \L_6, \L_{10}$ \\ 
\hline 
$2$ & $\L_4, \L_8, \L_{12}$ \\ 
\hline 
$3$ & $\L_2, \L_4, \L_9$ \\ 
\hline 
$4$ & $\L_3, \L_8, \L_{10}$ \\ 
\hline 
$5$ & $\L_3, \L_7, \L_9$ \\ 
\hline 
$6$ & $\L_3, \L_5, \L_9$ \\ 
\hline 
$7$ & $\L_3, \L_9, \L_{11}$ \\ 
\hline  
$8$ & $\L_1, \L_3, \L_9$ \\ 
\hline 
$9$ \& $10$ & $\L_5, \L_7, \L_{12}$ \\ 
\hline 
\end{tabular}
\end{minipage}
\hspace*{0.03\linewidth}
\begin{minipage}{0.3\linewidth}
\centering
\begin{tabular}{|r|p{3.5cm}|}
\hline
$j$ & Bundles intersecting at type $j$ vertices \\
\hline \hline
$11$ \& $12$ & $\L_1, \L_6, \L_{11}$ \\ 
\hline  
$13$ & $\L_2, \L_8, \L_{11}$ \\ 
\hline 
$14$ & $\L_1, \L_4, \L_{10}$ \\ 
\hline  
$15$ & $\L_2, \L_5, \L_8$ \\ 
\hline
$16$ & $\L_4, \L_7, \L_{10}$ \\ 
\hline 
$17$ & $\L_1, \L_5, \L_9$ \\ 
\hline 
$18$ & $\L_3, \L_7, \L_{11}$ \\ 
\hline 
$19$ & $\L_3, \L_5, \L_7, \L_9$ \\ 
\hline
$20$ & $\L_1, \L_3, \L_9, \L_{11}$ \\ 
\hline 
\end{tabular}
\end{minipage}
\hspace*{0.03\linewidth}
\begin{minipage}{0.3\linewidth}
\centering
\begin{tabular}{|r|p{3.5cm}|}
\hline
$j$ & Bundles intersecting at type $j$ vertices \\
\hline \hline
$21$ & $\L_3, \L_7, \L_9, \L_{11}$ \\ 
\hline   
$22$ & $\L_1, \L_3, \L_5, \L_9$ \\ 
\hline   
$23$ & $\L_1, \L_4, \L_7, \L_{10}$ \\ 
\hline   
$24$ & $\L_2, \L_5, \L_8, \L_{11}$ \\ 
\hline   
$25$ & $\L_1, \L_3, \L_7, \L_9, \L_{11}$ \\ 
\hline   
$26$ & $\L_1, \L_3, \L_5, \L_9, \L_{11}$ \\ 
\hline   
$27$ & $\L_3, \L_5, \L_7, \L_9, \L_{11}$ \\ 
\hline   
$28$ & $\L_1, \L_3, \L_5, \L_7, \L_9$ \\ 
\hline   
$29$ & $\L_1, \L_3, \L_5, \L_7, \L_9, \L_{11}$ \\ 
\hline   
\end{tabular}
\end{minipage}
\end{center}
\caption{Bundles intersecting at type $j$ vertices for $j=1,2,\ldots,29$.}
\label{tab:rec_L_j}
\end{table}

For $\lambda_6(m)$, all the $6$-wise crossings that are not 
at grid vertices, are at the centers of grid cells; we have 
\[ w_{29} = \area(\Gamma_1 \cap \Gamma_3 \cap \Gamma_5 \cap \Gamma_7 \cap \Gamma_9 \cap \Gamma_{11}) = 
\area(\Gamma_1 \cap \Gamma_5 \cap \Gamma_7 \cap \Gamma_{11}) = a_{12}. \]
It follows that
\[ \lambda_6(m) 
= \frac{a_6 + w_{29}}{\area(\sigma_0)} - O(m) 
= \frac{a_6 + a_{12}}{\area(\sigma_0)} - O(m) 
= \frac{2 m^2}{15} + \frac{m^2}{12} - O(m) = \frac{13 m^2}{60} - O(m). \]

Similarly for $\lambda_5(m)$, all the $5$-wise crossings 
that are not at grid vertices, \ie, types $25$ through $28$, 
are in the interiors of grid cells contained in eight small 
triangles inside $U$. For example, 
\[ w_{28} = \area(\Gamma_1 \cap \Gamma_3 \cap \Gamma_5 \cap \Gamma_7 \cap \Gamma_9 - \Gamma_{11})  
= \area(1,14,22) + \area(2,13,21) = 1/120. \]
Observe that sum of the areas of these eight small triangles 
equals to $a_{11}$. It follows that
\[ \lambda_5(m) 
= \frac{a_5 + \sum\nolimits_{j=25}^{28} w_j}{\area(\sigma_0)} - O(m) 
= \frac{a_5 + a_{11}}{\area(\sigma_0)} - O(m) 
= \frac{11 m^2}{30} + \frac{m^2}{30} - O(m) = \frac{2 m^2}{5} - O(m). \]

\begin{figure*}[htpb]
\centering
\includegraphics[scale=0.75]{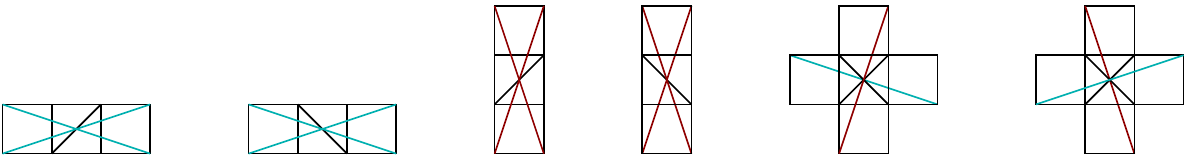}
\caption{These incidence patterns cannot occur.}
\label{fig:f19}
\end{figure*}

To estimate $\lambda_4(m)$, note that besides $4$-wise 
crossings at grid vertices, there are six types of $4$-wise 
crossings \ie, types $19$ through $24$, in the interiors of 
grid cells:
\begin{itemize} \itemsep 1pt
\item For types $19$ and $20$, there is one crossing per grid cell; and
\begin{align*}
w_{19} &= \area(\Gamma_3 \cap \Gamma_5 \cap \Gamma_7 \cap \Gamma_9 - \Gamma_1 - \Gamma_{11}) \\
&= (\area(2,10,13) - \area(2,10,21)) + (\area(9,14,22) - \area(1,14,22)) = 1/15,
\end{align*}
Type $20$ is a $90^{\circ}$ rotation of type $19$; therefore 
by symmetry, $w_{19} = w_{20}$.

\item For types $21$ and $22$, there is one crossing per grid cell; and
\begin{align*} 
w_{21} &= \area(\Gamma_3 \cap \Gamma_7 \cap \Gamma_9 \cap \Gamma_{11} - \Gamma_1 - \Gamma_5) \\
&= (\area(2,14,21) -\area(2,10,21)) + (\area(1,13,22) -\area(1,9,22)) = 1/40. 
\end{align*}
Type $22$ is the reflection in a vertical line of type $21$; 
therefore by symmetry, $w_{21} = w_{22}$.  

\item For types $23$ and $24$, there are four crossings 
per grid cell; and
\[ w_{23} = 4 \cdot \area(\Gamma_1 \cap \Gamma_4 \cap \Gamma_7 \cap \Gamma_{10}) 
= 4 \cdot \area(\Gamma_1 \cap \Gamma_7) = 4 \cdot \area(U_1) = 2/5. \]
Type $24$ is the reflection in a vertical line of type $23$; 
therefore by symmetry, $w_{23} = w_{24}$. 
\end{itemize}
Consequently, we have
\[ \lambda_4(m) = \frac{a_4 + \sum_{j=19}^{24} w_j}{\area(\sigma_0)} - O(m) 
= \Big(\frac{11}{15} + \frac{2}{15}+ \frac{1}{20} + \frac{4}{5} \Big)m^2 - O(m) 
= \frac{103 m^2}{60} - O(m). \]

Lastly, we estimate $\lambda_3(m)$. Besides $3$-wise crossings at 
grid vertices, there are $18$ types of $3$-wise crossings \ie, 
types $1$ through $18$, in the interior of grid cells:
\begin{itemize} \itemsep 1pt
\item For types $1$ and $2$, there is one crossing per 
grid cell; and
\[ w_1 = \area(\Gamma_2 \cap \Gamma_6 \cap \Gamma_{10}) 
= \area(\Gamma_2 \cap \Gamma_{10}) = \area(P(3,4,19,20)) = 1/4. \]
Type $2$ is a $90^{\circ}$ rotation of type $1$; therefore by 
symmetry, $w_1 = w_2$.

\item For types $3$ and $4$, there are two crossings per 
grid cell; and
\[ w_3 = 2 \cdot(\area(\Gamma_2 \cap \Gamma_4 \cap \Gamma_9)) 
= 2 \cdot(\area(P(3,4,7,8)) - \area(3,8,18) - \area(4,7,17)) = 1/2. \]
Type $4$ is the reflection in a vertical line of type $3$; 
therefore by symmetry, $w_3 = w_4$.

\item For types $5,6,7,8$, there is one crossing per grid cell; and
\[ w_5 = \area(\Gamma_3 \cap \Gamma_7 \cap \Gamma_9 - \Gamma_1 - \Gamma_5 - \Gamma_{11}) 
= \area(5,9,22) + \area(6,10,21)) = 1/20. \]
Type $6$ is the reflection in a vertical line of type $5$, 
and types $7$ and $8$ are $90^{\circ}$ rotations of types $6$ 
and $5$, respectively. Therefore by symmetry, 
$w_5 = w_6 = w_7 = w_8$. 

\item For types $9,10,11,12$, there is one crossing on 
the boundary of each grid cell; and
\[ w_9 = \area(\Gamma_5 \cap \Gamma_7 \cap \Gamma_{12}) 
= \area(\Gamma_5 \cap \Gamma_7) = \area(P(9,10,13,14)) = 1/6. \]
Type $10$ is the reflection in a horizontal line of type 
$9$, and types $11$ and $12$ are $90^{\circ}$ rotations of types $9$ 
and $10$, respectively. Therefore by symmetry, 
$w_9 = w_{10} = w_{11} = w_{12}$. 

\item For types $13,14,15,16$, there are four crossings per 
grid cell; and
\[ w_{13} 
= 4 \cdot (\area(\Gamma_2 \cap \Gamma_8 \cap \Gamma_{11} - \Gamma_5)) 
= 4 \cdot (\area(3,9,13) + \area(4,10,16)) = 1/5. \]
Type $14$ is the reflection in a vertical line of type $13$, 
and types $15$ and $16$ are $90^{\circ}$ rotations of types $13$ 
and $14$, respectively. Therefore by symmetry, 
$w_{13} = w_{14} = w_{15} = w_{16}$.

\item For types $17$ and $18$, there are two crossings per 
grid cell; and
\[ w_{17} = 2 \cdot (\area(\Gamma_1 \cap \Gamma_5 \cap \Gamma_9)) 
= 2 \cdot (\area(\Gamma_1 \cap \Gamma_5)) 
= 2 \cdot \area(P(1,2,9,10)) = 1/4. \]
Type $18$ is the reflection in a vertical line of type $17$; 
therefore by symmetry, $w_{17} = w_{18}$.
\end{itemize}

Consequently, we have
\begin{align*}
\lambda_3(m) &= \frac{a_3 + \sum_{j=1}^{18} w_j}{\area(\sigma_0)} - O(m) \\
& = \Big(\frac{4}{3} + \frac{1}{2} + 1 + \frac{1}{10} + \frac{1}{10}
+ \frac{1}{3} + \frac{1}{3}  + \frac{4}{5}  + \frac{1}{2} \Big)m^2 - O(m) = 5 m^2 - O(m). 
\end{align*}
\begin{table}[htpb]
\begin{center}
\begin{minipage}{0.24\linewidth}
\centering
\begin{tabular}{|l|r|}
\hline
$j$ & $w_j$ \\
\hline \hline
$1$ & $1/4$ \\ 
\hline 
$2$ & $1/4$ \\ 
\hline 
$3$ & $1/2$ \\ 
\hline 
$4$ & $1/2$ \\ 
\hline 
$5$ & $1/20$ \\ 
\hline 
$6$ & $1/20$ \\ 
\hline 
$7$ & $1/20$ \\ 
\hline  
\end{tabular}
\end{minipage}
\begin{minipage}{0.24\linewidth}
\centering
\begin{tabular}{|l|r|}
\hline
$j$ & $w_j$ \\
\hline \hline
$8$ & $1/20$ \\ 
\hline 
$9$ \& $10$ & $1/3$ \\ 
\hline
$11$ \& $12$ & $1/3$ \\ 
\hline  
$13$ & $1/5$ \\ 
\hline 
$14$ & $1/5$ \\ 
\hline 
$15$ & $1/5$ \\ 
\hline 
$16$ & $1/5$ \\ 
\hline  
\end{tabular}
\end{minipage}
\begin{minipage}{0.24\linewidth}
\centering
\begin{tabular}{|l|r|}
\hline
$j$ & $w_j$ \\
\hline \hline
$17$ & $1/4$ \\ 
\hline 
$18$ & $1/4$ \\ 
\hline 
$19$ & $1/15$ \\ 
\hline 
$20$ & $1/15$ \\ 
\hline
$21$ & $1/40$ \\ 
\hline 
$22$ & $1/40$ \\ 
\hline 
$23$ & $2/5$ \\ 
\hline  
\end{tabular}
\end{minipage}
\begin{minipage}{0.24\linewidth}
\centering
\begin{tabular}{|l|r|}
\hline
$j$ & $w_j$ \\
\hline \hline
$24$ & $2/5$ \\ 
\hline 
$25$ & $1/120$ \\ 
\hline 
$26$ & $1/120$ \\ 
\hline 
$27$ & $1/120$ \\ 
\hline 
$28$ & $1/120$ \\ 
\hline 
$29$ & $1/12$ \\ 
\hline
\end{tabular}
\vspace*{5mm}
\end{minipage}
\caption{Values of $w_j$ for $j=1, \ldots, 29$.}
\label{tab:rec_w_j}
\end{center}
\end{table}
The values of $\lambda_i(m)$, for $i=3,4,\ldots,12$, are summarized 
in Table~\ref{tab:lambda-12}; for convenience the linear terms are 
omitted. Since $m = n/12$, $\lambda_i$ can be also viewed as a function of $n$. 
\renewcommand*{\arraystretch}{1.3}
\begin{table}[H]
\begin{center}
\begin{tabular}{|c||c|c|c|c|c|c|c|c|c|c|}
\hline 
$i$ & $3$ & $4$ & $5$ & $6$ & $7$ & $8$ & $9$ & $10$ & $11$ & $12$ \\
\hline 
$\lambda_i(m)$ & $5 m^2$ & $\frac{103 m^2}{60}$ & $\frac{2 m^2}{5}$ 
& $\frac{13 m^2}{60}$ & $\frac{11 m^2}{105}$ & $\frac{13 m^2}{105}$
& $\frac{4 m^2}{105}$ & $\frac{m^2}{20}$ & $\frac{m^2}{30}$ 
& $\frac{m^2}{12}$ \\ 
\hline 
$\lambda_i(n)$ & $\frac{5 n^2}{144}$ & $\frac{103 n^2}{60 \cdot 144}$
& $\frac{2 n^2}{5 \cdot 144}$ & $\frac{13 n^2}{60 \cdot 144}$ 
& $\frac{11 n^2}{105 \cdot 144}$ & $\frac{13 n^2}{105 \cdot 144}$
& $\frac{4 n^2}{105 \cdot 144}$ & $\frac{n^2}{20 \cdot 144}$
& $\frac{n^2}{30 \cdot 144}$ & $\frac{n^2}{12 \cdot 144}$ \\
\hline
\end{tabular}
\caption{The asymptotic values of $\lambda_i(m)$ and $\lambda_i(n)$ for $i=3,4,\ldots,12$.}
\label{tab:lambda-12}
\end{center}
\end{table}
\renewcommand*{\arraystretch}{1}

The  multiplicative factor in Eq.\,\eqref{eq:recurrence-k}
is bounded from below as follows:   
\begin{align*} 
F(n) & \geq \prod_{i=3}^{12} B_i^{\lambda_i(n)} 
\geq 2^{\frac{5 n^2}{144}} 
\cdot 8^{\frac{103 n^2}{60 \cdot 144}} 
\cdot 62^{\frac{2 n^2}{5 \cdot 144}} 
\cdot 908^{\frac{13 n^2}{60 \cdot 144}} 
\cdot 24698^{\frac{11 n^2}{105 \cdot 144}} 
\cdot 1232944^{\frac{13 n^2}{105 \cdot 144}} \\
& \cdot 112018190^{\frac{4 n^2}{105 \cdot 144}} 
\cdot 18410581880^{\frac{n^2}{20 \cdot 144}} 
\cdot 5449192389984^{\frac{n^2}{30 \cdot 144}} 
\cdot 2894710651370536^{\frac{n^2}{12 \cdot 144}} 
\cdot 2^{-O(n)}. 
\end{align*}

We prove by induction on $n$ that $T(n) \geq 2^{cn^2 - O(n \log{n})}$ 
for a suitable constant $c>0$. It suffices to choose $c$ (using 
the values of $B_i$ for $i=3,\ldots,12$ in Table~\ref{tab:A_nB_n}) 
so that 
\begin{align*} 
\frac{1}{144} & \Big( 5 + \frac{103}{60} \log{8} + \frac{2}{5} \log{62} 
+ \frac{13}{60} \log{908} + \frac{11}{105} \log{24698} 
+ \frac{13}{105} \log{1232944} + \frac{4}{105} \log{112018190} \\ 
& + \frac{1}{20} \log{18410581880} + \frac{1}{30} \log{5449192389984} 
+ \frac{1}{12} \log{2894710651370536} \Big) 
\geq \frac{11c}{12}. 
\end{align*}

The above inequality holds if we set
\begin{equation} \label{eq:0.2053}
\begin{split}
c = \frac{1}{132} 
\Big( 5 & + \frac{103}{60} \log{8} + \frac{2}{5} \log{62} +
\frac{13}{60} \log{908} + \frac{11}{105} \log{24698} \\
& + \frac{13}{105} \log{1232944} + \frac{4}{105} \log{112018190}  + \frac{1}{20} \log{18410581880} \\ 
& + \frac{1}{30} \log{5449192389984} + \frac{1}{12} \log{2894710651370536} \Big) > 0.2053.
\end{split}
\end{equation}

\end{document}